\documentclass[11pt,reqno]{article}
\usepackage[utf8]{inputenc}
\usepackage{color}
\usepackage{mathrsfs}%                                          *
\usepackage{pifont}%                                            *
\usepackage{amsmath}%                                           *
\usepackage{amsthm}%                                            *
\usepackage{txfonts}%                                           *
 \usepackage{geometry}%                                         *
\usepackage{latexsym}%                                          *
\usepackage{amssymb}%                                           *
\usepackage{graphicx}%                                      
\usepackage{ifthen}%                                            *
\usepackage{fancyhdr}%                                          *
\usepackage{afterpage}%                                         *
\usepackage{xcolor}%                                            *
\usepackage{ifthen}%                                            *

\newlength{\defbaselineskip}
\newcommand{\setlinespacing}[2]%
           {\setlength{\baselineskip}{#1 \defbaselineskip}}
\newcommand{\doublespacing}{\setlength{\baselineskip}%
                           {1.5 \defbaselineskip}}

\newtheorem{theorem}{Theorem}

\newtheorem{cor}{Corollary}

\theoremstyle{definition}

\begin{document}
\doublespacing
\begin{center}
{\Large {\bf Sequential Detection of Transient Signals in High Dimensional Data Stream  }} \\
{\it Yanhong Wu{\footnote{Address correspondence to Yanhong Wu, Department of Mathematics, California State University Stanislaus, Turlock, CA 95382, USA; E-mail: ywu1@csustan.edu}} and D. Siegmund} \\
Department of Mathematics, California State University Stanislaus,\\
Turlock, California, USA \\
Department of Statistics, Stanford University, California, USA \\
\end{center}

\noindent{\bf ABSTRACT:} Motivated by sequential detection of transient signals in high dimensional data stream, we study the performance of EWMA, MA, CUSUM, and GLRT charts for detecting a transient signal in multivariate data streams in terms of the power of detection (POD) under the constraint of false detecting probability (FDP) at the stationary state. Approximations are given for FDP and POD. Comparisons show that the EWMA chart performs equally well as the GLRT chart when the signal strength is unknown, while its design is free of signal length and easy to update. In addition, the MEWMA chart with hard-threshold performs better when the signal only appears in a small portion of the channels. Dow Jones 30 industrial stock prices are used for illustration.  \\

\noindent{\bf Keywords:} CUSUM and EWMA procedures; False detection probability; MA and GLRT procedures; Power of detection; Transient signal.  \\
\noindent{\bf Subject Classifications:} 62L10; 60L15

\section{Introduction}
Consider a sequence of independent normal random variables $\{X_i\}$ that follow $N(0,1)$ for $i\leq \nu$ or $i> \nu+L$ and $N( \mu,1 )$ for $\nu +1 \leq i \leq \nu+L$ where the signal only appears on the interval $[\nu +1, \nu+L]$, called transient signal with length $L$. The traditional sequential change detection problem with $\mu>0$ assumes $L=\infty$. the Shewhart chart (Shewhart (1931)) makes an alarm when $X_i>c$ (or $(X_{(i-1)k+1}+...X_{ik})/\sqrt{k} >c$ for group size $k$).  That means, the alarm is only raised based on the current (grouped) observation. The CUSUM (Page (1954)), EWMA (Roberts (1959)), and Shiryayev-Roberts (S-R) (Roberts (1966) and Shiryayev (1963)) charts are developed for a quick detection of change. For example, we may focus to find a stopping time $\tau$ that minimizes certain criterion based on the conditional delay detection time $CADT(\nu) = E[\tau-\nu|\tau >\nu]$ for given $ARL_0 =T$, where the optimal properties of the CUSUM and S-R charts have been established under different criteria for a given reference value $\delta $ of $\mu$.  
When $L=1$, i.e. the signal just appears in one observation, the Shewhart chart will possess certain optimal property as shown in Pollak and Krieger (2013) and Moutstakides (2014). 

 The sequential detection of transient signal requires new criterion to evaluate a monitoring chart. Both the signal strength and signal length are the concerns in addition to the change time.  An obvious chart is the MOSUM (moving sum) chart (Bauer and Hackel (1978)) or the MA (moving average) chart that uses average of the partial sum $\sum_{n-w+1}^{n} X_i /w =\bar{X}_{n;w}$ to detect the signal by raising an alarm at
\[
\tau_{MA} =\inf\{n>0: \bar{X}_{n;w} > h \}  
\]
where $w$ is a pre-selected window size that can be treated as a reference value for the signal length $L$ and $h$ is the reference value for the signal strength. Several authors have proposed other procedures. General weighted moving average chart is considered in Lai (1974).  Gueppie et al. (2012) proposed to use a window-limited CUSUM procedure under the constraint on the false alarm probability $P_0(\nu <\tau<\nu+L)$ for just one alarm interval where $P_0(.)$ is the probability without change. Noonan and Zhigljavsky (2020) studied the power of the MA chart by using the corrected diffusion approximation under the constraint on $ARL_0$. Tartakovsky et al. (2021) considered a modified CUSUM procedure that is shown to be optimal in terms of average delay detection time when the signal length $L$ follows an exponential distribution under the constraint of the conditional false alarm probability $\max_{\nu} P_0[\tau \leq \nu +L |\tau>\nu ]$. 

To evaluate the performance of a monitoring chart for transient signals, we consider the following two scenarios. First, we can assume that data from consecutive periodic intervals of length $T$ $((k-1)T, kT]$ for $k=1,2,...$ are monitored subject to a constraint on $P_0(\tau \leq T) \leq \alpha _0$ and expect signal appears randomly in one interval. Second, we can fix $ARL_0 =E_0 \tau $ as a measure of false alarm rate in the long run and observe consecutive false alarm cycles until the signal appears in one alarm cycle. Here $1/ARL_0$ is the false alarm rate in the long run. In both situations, the power of detection $P_{\mu}^{\nu}[\tau \leq \nu +L |\tau > \nu]$ is to be evaluated. For the change point detection, the performance is evaluated based on the stationary average delay detection time (SADDT) as considered in Srivastava and Wu (1993),  Reynolds and Stoumbos (2004), and Knoth (2021) for numerical and theoretical evaluations. 

To unify our discussion, we shall assume that the monitoring  process will approach stationary. Thus the false alarm time is approximately exponential. This implies $P_0(\tau \leq T) \approx 1-\exp(-T/ARL_0)$. That means, fixing $ARL_0$ is equivalent to fixing $P_0(\tau \leq T)$ for any $T$. On the other hand, at the change point $\nu$, we can assume the controlled monitoring process is at the stationary state. 
Denote $P_{\mu} ^*(.)$ as the probability measure when the controlled process starts at the stationary state with signal strength $\mu$. That means when $\nu \rightarrow \infty$, $P_{\mu}^{\nu}[\tau \leq \nu +L |\tau > \nu] \approx P_{\mu}^*[\tau \leq L] $. By noting that for any $T$,
\[ \alpha (T) = P_0^*(\tau \leq T) =1-P_0^*(\tau >T) = 1- (P_0^*(\tau>L))^{T/L}=1- (1-\alpha(L))^{T/L}, \] 

we can evaluate and compare a monitoring procedure based on the {\em false detection probability} (FDP) defined as \[\alpha (L)= P_0^*(\tau \leq L),\]  
and the {\em power of detection} (POD) defined as \[\beta(L, \mu)=  P_{\mu}^*(\tau \leq L)\]
for signal length $L$.  Obviously, $\alpha (L)= \beta(L, 0)$.  Just like in standard statistical hypothesis test, we can evaluate the power of detection under the constraint on false detection probability. In addition, we can also study the average delay detection time $E_{\mu}^*[\tau|\tau\leq L]$, conditioning on the signal being detected.

In this communication, we evaluate and compare several monitoring procedures based on the EWMA, moving average, CUSUM, and generalized likelihood ratio tests with the transient signal. We show that the localization method developed in Siegmund and Yakir (2000), Siegmund, et al. (2010, 2011), and Xie and Siegmund (2013) can be naturally extended here and gives accurate approximations for FDP. It also justifies for simpler approximations by using boundary correction technique.  In addition, we consider the situation when the signal only appears in a portion of channels. In Section 2, we first present the results for the one-dimensional case based on EWMA, MA, CUSUM, and GLRT charts. For fixed FDP, comparisons of POD and average delay detection time given the signal being detected are conducted. 
In Section 3, we first present the multivariate generalizations for the EWMA, MA,  CUSUM, and GLRT procedures that are signal direction invariant. Then we give approximations for FDP. Comparisons of PODs are conducted. In Section 4, we consider modified multivariate monitoring charts with adaptive weights  in signal strength in order to increase POD when signal only appears in a small portion of the channels. Our main finding is that when the signal length $L$ is unknown, EWMA chart (with hard threshold) should be recommended due to its simple calculation and overall better performance in terms of POD, particularly when signal only appears in a portion of the channels. Dow Jones 30 industrial stock prices are used for illustration.  

\section{One-dimensional Monitoring Charts}
In this section, we present the results for several monitoring charts for one-dimensional observations. 
\subsection{Definitions}

We shall continue using the notations in the introduction section. 

\noindent{\bf (1) EWMA chart:} Define $Y_n =(1-\beta)Y_{n-1} +\beta X_n$ as the EWMA process for $Y_0=0$. The EWMA chart makes an alarm at
\[
\tau_{EW}=\inf \{ n>0:  Y_n >b (\frac{\beta}{2-\beta})^{1/2} \}, 
\]
where $ \frac{\beta}{2-\beta} =\lim_{n\rightarrow \infty} Var (Y_n) $.

\noindent{\bf (2) MA chart:}  For window length $w$, define
\[
\tau_{MA} =\inf\{n \geq w: \bar{X}_{n;w} =(X_{n-w+1}+...+X_n)/w > h \}.
\]

\noindent {\bf Remark:} To increase the power of the detection for the EWMA chart for transient signal, we may consider the following moving-EWMA chart.

 Define $Y_{n;w} =\sum_{k=1}^w \beta (1-\beta)^{k-1} X_{n-k+1}$. The moving-EWMA chart makes an alarm at
\[
\tau_{MEW}=\inf \{ n\geq w:  Y_{n;w}/(1-(1-\beta)^w) >h \}, 
\]
where $ 1-(1-\beta)^w $ is the total weight. Note that moving-EWMA is a special case of generalized weighted moving average (GWMA) considered in Lai (1974). 

\noindent{\bf (3) CUSUM procedure:} For the reference signal strength $\delta $, the CUSUM procedure makes an alarm at
\[
\tau_{CS} =\inf\{ n>0:  Y_n =\max(0, Y_{n-1} +X_n-\delta/2) >d \},
\]
where $Y_0 =0$. 

For unknown signal strength, Siegmund and Venkatraman (1995) considered the generalized likelihood ratio (GLR) statistic under the change point model that makes an alarm at
\[
\tau_{GLR}= \inf\{ n>0: \max_{0\leq w<n} \sqrt{w}\bar{X}_{n;w}>b\}. 
\]

For finite signal length, we may consider an interval $w_0 \leq L \leq w_1$ as possible candidates for signal length. As Lai (1995) proposed, we can consider the following windowed GLR procedure that can also be considered as a generalization of the MA chart. 

\noindent{\bf (4) windowed-GLR chart:} For a possible signal window $[w_0, w_1]$ for the signal length, an alarm is made at
\[
\tau_{WGL}= \inf\{ n>w_1: \max_{w_0\leq w<w_1} \sqrt{w}\bar{X}_{n;w}>b\}. 
\]

\subsection{Approximations for FDP}

In the following, we evaluate the FDP for the above four monitoring charts. 
Note that for EWMA and MA charts, $b (\frac{\beta}{2-\beta})^{1/2}$ and $h$ play the role as the minimum signal strength we want to detect and can be assumed quite small. The following theorem gives the approximations for the FDP of EWMA and MA charts and the proofs can be obtained by slightly modifying the ones given in Siegmund, et al. (2011) and are given in the appendix. 

\begin{theorem}
(i) For the EWMA procedure, as $\beta \rightarrow 0$, $\beta L \rightarrow \infty$, and $b \rightarrow \infty$ such that $b^2 \beta \rightarrow $ constant, 
\begin{equation}
\alpha = P_0^*[\max_{1\leq k \leq L} Y_k \geq b (\beta /(2-\beta))^{1/2}] 
\approx L\beta b \phi (b) \upsilon( b\sqrt{2\beta})
\approx L\beta b \phi (b) e^{-b\rho_+ \sqrt{2\beta}} .
\end{equation}
(ii) For the MA procedure, as $w\rightarrow \infty$ and $L \rightarrow \infty$,
\begin{equation}
\alpha = P_0^*[ \max_{1\leq k \leq L} \frac{X_{k-w+1}+...+X_k}{w} \geq h] 
\approx 
\frac{Lh}{\sqrt{w}} \phi(h\sqrt{w}) \upsilon(\sqrt{2}h) 
\approx \frac{Lh}{\sqrt{w}}  \phi(h\sqrt{w}) e^{-\sqrt{2} h\rho+}. 
\end{equation}
where $\upsilon(x) \approx e^{-\rho_+ x} $ is given in Siegmund (1985, Pg.82) with $\rho_+ \approx 0.5826$ and $\phi(x)$ is the standard normal density function.  
\end{theorem}

The above approximations can be used to obtain a simple boundary correction to match the results in the continuous time case. For (i), $b$ is corrected to $b^* =b+ \rho_+ \sqrt{2\beta} $ and for (ii), $h$ is corrected to $h^* =h +\sqrt{2} \rho_+ /w$ (see Wu and Wu (2022)). For the GLR chart, Siegmund and Venkatraman (1995) obtained the following approximation:
\begin{theorem} As $b\rightarrow \infty$ and $m\rightarrow \infty$ such that $m=cb^2$ for some $0<c<\infty$, then 
\[
P_0[\tau_{GLR} \leq m] \approx mb\phi (b) \left (\int_{c^{-1/2}}^{\infty} x \upsilon^2 (x) dx -c^{-1} \int_{c^{-1/2}}^{\infty} x^{-1} \upsilon^2(x) dx \right ). 
\]
\end{theorem}

The following theorem gives the approximation for the GWMA chart or the windowed-GLR chart and the proof is based on the technique used in Siegmund, et al. (2010). 

\begin{theorem}
As $b\rightarrow \infty$ and $L \rightarrow \infty$ such that $\sqrt{w_0}/b, \sqrt{w_1}/b $ go to constants, 
\begin{equation}
P_0^*(\tau_{WGL} \leq L) \approx b^3 \phi(b) \int_{W_0/L}^{W_1/L} \upsilon^2(\frac{b}{\sqrt{Lt}}) \frac{1}{{4t^2}} dt 
=L b\phi(b) \int_{b/\sqrt{W_1}}^{b/\sqrt{W_0}} \frac{u \upsilon^2 (u)}{2} du.
\end{equation}
\end{theorem}

\noindent{\bf Proof}: We only briefly mention the main steps. First, define 
\[
Z_t^w =\frac{1}{\sqrt{w}} \sum_{t-w+1}^t X_i.
\]
We define a changed measure $P_t^w(.)$ such that under $P_t^w(.)$, $X_{t-w+1} ,..., X_t$ are i.i.d following $N(\theta /\sqrt{w}, 1)$ and other $X_i$ follow N(0,1) for other $i$'s. Then the likelihood ratio of $X_1,..., X_L$ under $P_t^w(.)$ with respect to $P_0^*(.)$ equals to 
\[
\frac{dP_t^w}{dP_0^*} (Z_t^w) =exp(\theta Z_t^w -\theta^2/2).
\]
Let $\theta = b$, $\tilde{l_t^w} = \theta (Z_t^w- \theta )$, and $l_t^w =\theta Z_t^w-\theta^2/2$. 
Then the localization technique under the changed measure gives the following approximation
\[
P_0^*(\tau_{WGL} \leq L ) \approx \frac{1}{b}\phi(b)
\sum_{1\leq t \leq L} \sum_{w_0\leq w \leq w_1} E_t^w[ \frac{M_t^w}{S_t^w}],
\]
where $M_t^w =\max_{1\leq s \leq L; w_0\leq u\leq w_1} \exp(l_s^u-l_t^w))$ and $S_t =\sum_{1\leq s \leq L; w_0\leq u\leq w_1} \exp(l_s^u-l_t^w)$. 

Next, when $w$ is large and $u$ is not too far away from $w$, we can write for $s<t$
\begin{eqnarray*}
\theta(Z_s^u-Z_t^w )&= & \theta( \frac{1}{\sqrt{u}} (\sum_{s-u+1}^s X_i -\sum_{t-w+1}^t X_i ) +Z_t^w (\sqrt{w/u} -1)) \\
&\approx & \theta( \frac{1}{\sqrt{u}} ( \sum_{s-u+1}^{t-w}X_i - \sum_{s+1}^t X_i ) +Z_t^w \frac{w-u}{2u} ) \\
& = & \frac{\theta}{\sqrt{u}} (\sum_{s-u+1}^{t-w} (X_i -\frac{Z_t^w}{2\sqrt{u}} ) 
-\frac{\theta}{\sqrt{u}}  \sum_{s+1}^t (X_i - \frac{Z_t^w}{2\sqrt{u}} )) .
\end{eqnarray*}
For  the first summation, it is a random walk with variance $ \theta^2/ u$ from the first term and drift $-\theta^2/2u$ from the second term. For the second summation, it has variance $\theta^2/u$ and drift $-\theta^2/2u$ as $X_i $ by noting that $Z_t^w$ has conditional mean $Z_t^w/\sqrt{w}$ and  $w$ is not too far away from $u$. As the covariance between the two random walks is at the lower order, so it can be treated as two independent random walks with different time scales. Similar argument applies when $s>t$. From Siegmund, et al. (2010), 
\[
E_t^w[ \frac{M_t^w}{S_t^w}] \approx  (\theta^2/2w)^2 \upsilon^2 (\theta/\sqrt{w} ) = (b^4/4w^2) 
\upsilon^2 (b/\sqrt{w}).
\]
The final approximation is obtained by changing the summation as an integration. \qed

For the CUSUM procedure, Siegmund (1988) derived a second order approximation when $Y_0=0$.

As $b\rightarrow \infty$ and $L \rightarrow \infty$ such that
$L(\delta/2)/d \geq 1+\epsilon$ and $L^2 exp(-\delta d) \rightarrow 0$,
\[
P_0(\tau_{CS} \leq L|Y_0=0) \approx (\delta (\delta L/2 -(d+2\rho_+))+3)e^{-\delta (d+2\rho_+)} .
\]

Since the CUSUM process is a renewal process at the points when it goes back to zero, the only modification is on the first interval where the starting point can be treated following the stationary distribution: $P_0^*(Y_0 =0) \approx 1-exp(-\rho_+ \delta )$ and
\[
P_0^*(Y_0 >x|Y_0>0) \approx  exp(-\delta x ) .
\]
Instead of modifying the above approximation, we can directly use the approximation for $ARL_0$ (Siegmund (1985, pg. 27)) to estimate the rate for the exponential distribution. For $\alpha =0.01$ and $L=20$, $ARL_0 \approx 2000$. The POD can be simulated with the above approximated stationary distribution.

Here, a similar approach following the lines in Theorems 1 and 3 gives the following result and will be used for multi-dimensional case.
\begin{theorem}
{\it As $d\rightarrow \infty$ and $L \rightarrow \infty$ such that and $L \exp(-\delta d) \rightarrow 0$,
\[
P_0^*(\tau_{CS} \leq L) \approx (L\delta^2 /2)e^{-\delta (d+2\rho_+)} .
\]
}
\end{theorem}

\noindent{\bf Proof.} 
Note that $\tau_{CS}$ can be written as
\[
\tau_{CS} =\inf\{t>0: \max_{0\leq w \leq t} w (\bar{X}_{t,w} -\delta/2) >d\}
\]
Define the changed measure $P_t^w (.)$ such that under $P_t^w (.)$ $X_{t-w+1}, ..., X_t$ follow $N(\delta, 1)$ and other $X_i$'s follow $N(0,1)$. Then
\[
\frac{P_t^w}{dP_0^*} ( \bar{X}_{t,w} )= \exp( w\delta (\bar{X}_{t,w} -\delta/2)).
\]
Denote by $l_t^w =\delta w(\bar{X}_{t,w} -\delta/2)) $. We can get
\[
P_0^*(\tau_{CS} \leq L) \approx e^{-\delta d} 
\sum_{1\leq t \leq L} \sum_{0\leq w \leq t} E_t^w[ \frac{M_t}{S_t} 
e^{-(l_t^w -\delta d) +log(M_t)}; l_t^w +\ln M_T >\delta d ] ,
\]
where $M_t =\max_{1\leq s \leq L} \max_{0\leq u \leq t}(l_s^u -l_t^w ) $ and 
$ S_t =\sum_{1\leq s \leq L} \sum_{0\leq u \leq t}(l_s^u -l_t^w )  $. 

Under $P_t^w(.)$, $\delta (\bar{X}_{t,w} -\delta/2)$ has mean $\delta^2/2$. Denote by $w_0= 2d/\delta$, then for any $ w$ such that $w =2d/\delta +2 \mu \sqrt{w_0} /\delta^2$, $(l_t^w -\delta d)/\sqrt{w} \rightarrow N(\mu, \delta^2)$. 
From Siegmund, et al. (2010), 
\[
\sum_{0\leq w \leq t} E_t^w[ 
e^{-(l_t^w -\delta d) +log(M_t)}; l_t^w +\ln M_T >\delta d ] \approx \int \frac{1}{\sqrt{w}{d}} \phi (\frac{\delta^2 w/2 -\delta d}{\sqrt{w}d } ) d w =2/\delta^2.
\]
On the other hand, $l_s^u -l_t^w$ can be approximated as the sum of two independent normal random walk with drift $ \delta^2 /2$ and variance $\delta^2$. So $E_t^w[ \frac{M_t}{S_t}] \approx (\delta^2/2)^2 \upsilon^2( \delta )$.  Thus, 
\[
P_0^*(\tau_{CS} \leq L) \approx L e^{-\delta d} \frac{\delta^2}{2}\upsilon^2(\delta). \qed
\]

The simulation results show that the approximations for the EWMA and MA charts are quite good, while they underestimate the true values for the CUSUM and GLRT charts. 

\subsection{POD and comparison}

In Table 1, we compare the PODs between EWMA, MA, CUSUM, and GLRT charts for the one-sided change for $L=20, 30, 40, 50$ with $\alpha =0.01,0.015,0.02, 0.025$, respectively.  All simulations are replicated 50,000 times. $L=20$ is taken as the benchmark for design with FDP = 0.01. \\
(i) For the EWMA chart, we take the commonly used weight $\beta =0.05$ and alarm limit $b=2.95$. \\
(ii) For the MA chart, we take $w=10, 20, 50$ that correspond to $h=0.99074$, $0.6578$ and 0.394 by using the approximation with $FDP =0.01$.  \\
(iii) For the CUSUM procedure, we use $b=10.8$ for $\delta =0.5$ and $b=5.88$ for $\delta =1.0$.\\
(iv) For the GLRT chart, we take $b=3.27$ for $20<w\leq 50$.

The numbers in the bracket are the approximations by using Theorems 1, 3, and 4. We see that the MA chart heavily depends on the window width. As expected, the MA chart performs better when $L=w$. Overall, the EWMA chart performs as well as the GLRT and both should be recommended for practical use. However, the EWMA chart is easier to calculate and does need to known the signal length. 

\begin{table}[ht]
\begin{center}
Table 1. Comparison of POD for one-sided change\\
\begin{tabular}{|c| c |c |c | c|c | c|c | c| }   \hline
L & $\mu$ & EWMA & MA(10) & MA ($20$) & MA(50) &CUSUM ($0.5$) & CUSUM ($1.0$) &GLRT\\ \hline
20 &0.0 & 0.0105 & 0.0090& 0.0105 & 0.0102 &0.0096 & 0.0106& 0.00984\\
  & & (0.0103) & (0.0082) & (0.0090) &(0.0066) & (0.0063) & (0.0087) & (0.0049)  \\ 
& 0.25 & 0.0634 &0.0576 & 0.0719 & 0.0329 &0.0516 & 0.0619 &0.0507  \\
& 0.5 & 0.2641 & 0.2387 & 0.3188 & 0.1088 & 0.2363 & 0.2742&0.2401 \\
& 0.75& 0.6232 & 0.5791 & 0.7180 & 0.2791& 0.6123 & 0.6503& 0.6167\\
& 1.00 & 0.9043 & 0.8750 & 0.9516 &0.5380&  0.9076 & 0.9214 & 0.9081\\
& 1.25 & 0.9892 & 0.9842 & 0.9970 &0.7870& 0.9911 & 0.9932& 0.9910 \\
& 1.50 & 0.9994& 0.9994 & 0.99992 &0.9300 & 0.99962 & 0.99974& 0.99976 \\
& 1.75 & 0.99994 & 0.999998 & 1.0000 & 0.9853& 1.0000 & 1.0000& 1.00\\
& 2.0 & 1.0 & 1.0 & 1.0 & 0.9979 & 1.0 & 1.0 &1.00  \\ \hline
30& 0.0 & 0.0143 & 0.0125 & 0.0146 &0.0138& 0.0133 & 0.0156 &0.0142 \\
 &0.25 & 0.1228 &0.0908 & 0.1275 &0.0748&  0.1068 & 0.1016& 0.1153\\
& 0.5 & 0.4998 & 0.3637 & 0.5103 &0.03065&  0.4779 & 0.4322& 0.5078 \\
& 0.75& 0.8889 & 0.7612 & 0.8866 &0.6935& 0.8876 & 0.8430& 0.9007 \\
&1.00 & 0.9926 & 0.9668 & 0.9935 & 0.9382& 0.9937 & 0.9878 & 0.9954\\
&1.25 & 0.9999 & 0.9985 & 0.99992 &0.9948&  0.99996 & 0.99978 & 0.9999 \\
& 1.50 & 1.0000 & 0.99998 & 1.0000 &0.9999& 1.0000 & 1.0000 & 1.00\\
&1.75 & 1.0000& 1.0000 & 1.0000 &1.0000& 1.0000 & 1.0000 &1.00 \\
& 2.0 & 1.0 & 1.0 & 1.0 & 1.0 &  1.0 & 1.0 &1.00   \\ \hline
40 &0.0 & 0.0169 & 0.01754 & 0.01998 &0.01744 & 0.0169 & 0.0190 &0.0165 \\
 &0.25 & 0.1882 &0.1243& 0.1814 &0.1394&  0.1676 & 0.1406& 0.1764\\
 &0.5 & 0.6833& 0.4680 & 0.6420 &0.5924& 0.6750 & 0.5518& 0.6867\\
 &0.75& 0.9702 & 0.8625 & 0.9563 &0.9424&  0.9715 & 0.9308& 0.9745 \\
 &1.00 & 0.9996 & 0.9907 & 0.9990 &0.9983& 0.9997 & 0.9983& 0.9997 \\
 &1.25 & 1.0000 & 0.99996 & 1.0000 &1.0000& 0.99998 & 0.99998& 1.0000\\
 &1.50 & 1.0000 & 1.0000 & 1.0000 &1.0000& 1.0000 & 1.0000 & 1.00 \\
 &1.75 & 1.0000& 1.0000 & 1.0000 & 1.0000& 1.0000 & 1.0000 &1.00\\
 & 2.0 & 1.0 & 1.0 & 1.0 & 1.0 & 1.0 & 1.0 &1.00 \\ \hline
50 & 0.0 & 0.0217& 0.0217 & 0.0253 &0.2090 &0.0203 & 0.0233& 0.0202 \\
 &0.25 & 0.2591 &0.1550 & 0.2309 &0.2447&  0.2353 & 0.1802 &0.2356 \\
 &0.5 & 0.8093 & 0.5549 & 0.7461 &0.8351 & 0.8020 & 0.6540 & 0.8051 \\
 &0.75& 0.9923 & 0.9202 & 0.9839 &0.9962 & 0.9940 & 0.9694 & 0.9931\\
 &1.00 & 0.99998 & 0.9975 & 0.99998 &1.0000& 1.0000 & 0.99972& 1.0000 \\
 &1.25 & 1.0000& 1.0000 & 1.0000 &1.0000 &1.0000 & 1.0000& 1.0000 \\
 &1.50 & 1.0000 & 1.0000 & 1.0000 & 1.0000& 1.0000 & 1.0000 &1.0000 \\
 & 1.75 & 1.0000& 1.0000 & 1.0000 & 1.0000& 1.0000 & 1.0000 & 1.0000\\
 & 2.0 & 1.0 & 1.0 & 1.0 & 1.0 & 1.0& 1.0 & 1.0 \\ \hline
\end{tabular}
\end{center}
\end{table}

In the following, we derive some asymptotic results for POD. 

(1) For the EWMA chart, under $P_{\delta}^*(.)$ we can write
\[
\beta(L,\delta )= P_{\delta}^*(\tau){EW} \leq L)= P_0^*(\max_{0\leq t \leq L} (Y_t + (1-(1-\beta)^t)\delta) >b\sqrt{\beta/(2-\beta)} ).
\]
To obtain a theoretical expression for the POD, we can take $\beta$ as the time unit and assume $\beta \rightarrow 0 $, $b\rightarrow \infty $, and $L\rightarrow \infty $ such that $b\sqrt{\beta/(2-\beta)} \rightarrow h$ and $L\beta \rightarrow \infty$.

We consider two situations. In the local case when the signal strength $\delta <h$, we can use the approximation of FDP for an interval around $\beta L =t $. That means, let
\[
b_t =b-(1-e^{-t})\delta/\sqrt{\beta/(2-\beta)} .
\]
Thus,
\begin{equation}
\beta(L,\delta ) \approx \int_0^{\beta L} \beta b_t \phi(b_t) \upsilon(b_t \sqrt{2\beta}) dt \approx 1-\exp(-\int_0^{\beta L} \beta b_t \phi(b_t) \upsilon(b_t \sqrt{2\beta}) dt ).
\end{equation}
See Wu and Wu (2022) for a detailed proof and some numeric comparisons. The second is when the signal strength $\delta \geq h$. By ignoring the overshoot at the crossing time, we have
\[
Y_{\tau_{EW}} + (1-e^{-\beta \tau_{EW}})\delta= h.
\]
This gives
\[
\beta\tau_{EW} =- \ln (1-\frac{h}{\delta} +\frac{1}{\delta} Y_{\tau_{EW}} ) .
\]
Note that the boundary crossing is mainly caused by the signal and $Y_{\tau_{EW}} =\sqrt{\beta/(2-\beta)} Z_{\tau_{EW}}$. A second order expansion gives
\[
\beta\tau_{EW} = -\ln (1-\frac{h}{\delta}) -\frac{Y_{\tau_{EW}}}{\delta-h} +\frac{Y_{\tau_{EW}}^2}{2(\delta -h)^2} +o_p(\beta).
\]
Thus, 
\begin{eqnarray*}
E_{\delta}^*(\tau_{EW}) &=& -\frac{1}{\beta} \ln (1-\frac{h}{\delta}) +\frac{1}{4(\delta -h)^2} +o(1 ). \\
Var_{\delta}^*(\tau_{EW} ) &=& \frac{1}{2\beta (\delta-h)^2} +O(1). 
\end{eqnarray*}
The normal approximation can be used to approximate the POD as
\begin{equation}
\beta(L,\delta) \approx \Phi \left (\frac{\beta L+\ln(1-h/\delta) -\beta /4(\delta-h)^2}{\sqrt{\beta/2}/(\delta -h)} \right ) .
\end{equation}
A more accurate approximation can be obtained by using the continuous correction boundary for $b $ with $b^* =b+ \sqrt{2\beta} \rho_+$ or $h^*=h+\beta \rho_+ $ in the approximation for FDP in Theorem 1. 

(2) Next, we consider the MA chart. Here $1/w$ is treated as the time scale. In the local case, we assume $w\rightarrow \infty$ and $\delta <h$. The POD can be approximated as
\begin{eqnarray}
P_{\delta}^*(\tau_{MA} \leq L ) = && P_0^*( \max_{1\leq t \leq L} (\bar{X}_{t;w} + \min(t/w,1) \delta ) >h ) \nonumber \\
= &&    \int_0^{L/w} (\sqrt{w} (h - \min(u, 1) \delta) \phi( (\sqrt{w} (h - \min(u, 1) \delta) e^{-\sqrt{2} \rho_+  (h - \min(u, 1) \delta)} du  \\
 = &&   \int_0^1 (\sqrt{w} (h - u \delta) \phi( (\sqrt{w} (h - u\delta)) e^{-\sqrt{2} \rho_+ (h - u\delta)}  d u \nonumber\\
  & &     + (L/w-1)( \sqrt{w}(h - \delta) \phi (\sqrt{w}(h - \delta))) e^{-\sqrt{2} \rho_+ (h - \delta)}, ~ if ~L\geq w. \nonumber
\end{eqnarray}

When $\delta \geq h$, by ignoring the overshoot, at the crossing time
\[
\bar{X}_{t;w} + \min(t/w,1) \delta =h , 
\]
i.e. $
t/w= h/\delta - \bar{X}_{t;w}/\delta .$ Therefore, as $w\rightarrow \infty$, 
\[
E_{\delta}^*[\tau_{MA}] =h w /\delta, ~~~Var_{\delta}^* (\tau_{MA}) = w/\delta^2. 
\]
The POD can be approximated as
\begin{equation}
P_{\delta}^*(\tau_{MA} \leq L ) \approx \Phi ( \delta \sqrt{w} (L/w -h/\delta )). 
\end{equation}
Similar to the EWMA chart, we can use a corrected boundary $h^*= h+\sqrt{2}\rho_+/w$. 

(3) For the CUSUM chart, note that the stopping time $\tau_{CS}$ can be written as 
\[
\tau_{CS} =\inf\{t>0: \max_{0\leq w \leq t} w(\bar{X}_{t;w} -\delta/2) >d\}.
\]
Suppose the signal $\mu >\delta/2$ appears at the change time $\nu$ for $\nu <t$. By ignoring the overshoot, at the crossing time,  we have
\[
\max_{0\leq w \leq t} w(Z_{t;w}/\sqrt{w} +\mu \min(t-\nu/w,1) -\delta/2) = d.
\]
As $d\rightarrow \infty$, the optimal $w \rightarrow \infty$. Therefore at the first order,
$w =t-\nu$. So
\[
(t-\nu)(\mu-\delta/2) +\sqrt{t-\nu} Z_{t;w} =d.
\]
This gives $
t-\nu =\frac{d}{\mu-\delta/2} (1- \frac{\sqrt{t-\nu}}{d} Z_{t;w} ) $. Thus, 
\begin{eqnarray*}
\sqrt{t-\nu} &=& (\frac{d}{\mu-\delta/2})^{1/2} (1- \frac{\sqrt{t-\nu}}{2d} Z_{t;w} +o_p(1/d)) \\
&= & (\frac{d}{\mu-\delta/2})^{1/2} -\frac{1}{2(\mu-\delta/2)} Z_{t;w} +o_p(1/\sqrt{d} )
\end{eqnarray*}
So
\[
t-\nu =\frac{d}{\mu-\delta/2} -\frac{\sqrt{d}}{(\mu-\delta/2)^{3/2}} Z_{t;w} +
\frac{Z_{t;w}^2}{2(\mu-\delta/2)^2}+o_p(1).
\]
This gives
\begin{eqnarray*}
E_{\mu}^*[\tau_{CS}] &=& \frac{d}{\mu-\delta/2}+\frac{1}{2(\mu-\delta/2)^2}+o(1), \\
Var_{\mu}^*(\tau_{CS}) &= &d/(\mu-\delta/2)^3 +O(1).
\end{eqnarray*}
An approximate for the POD can be obtained as 
\[
P_{\mu}^*(\tau_{CS}\leq L) \approx \Phi \left(\frac{L-d/(\mu-\delta/2) -1/(2(\mu-\delta/2)^2)}{\sqrt{d}/(\mu-\delta/2)^{3/2}} \right).
\]

(4) The derivation for GLRT chart is relatively easy. By noting at the alarm time,
\[
\max \sqrt{w} (Z_{t;w}/\sqrt{w} +\mu \min ((t-\nu)/w,1)) =b,
\]
we have $t-\nu= w$ and 
\[
t-\nu =(b-Z_{t;w})^2/\mu^2 = b^2/\mu^2 -2bZ_{t;w}/\mu^2 +Z_{t;w}^2/\mu^2.
\]
Thus, if $b/\sqrt{w_1} < \mu <b/\sqrt{w_0}$, 
\[
E_{\mu}^*(\tau_{GL}) =(b^2+1)/\mu^2; ~~~ Var_{\mu}^* (\tau_{GL}) = 4b^2/\mu^4 .
\]
An approximation can be obtained as
\[
P_{\mu}^*(\tau_{GL} \leq L) \approx \Phi ((L-(b^2+1)/\mu)/(2b/\mu^2 )).
\]

\section{Multivariate Monitoring Charts}
\subsection{Definitions}
We consider the following multivariate transient change model with possible cross-dependence. 
Let \[
X_{it} = \mu_{it} +Z_{it} ,
\]
for $i=1,...,N$ and $t=1,2,...$, where $\mu_{it}=\mu_i I_{[\nu \leq t \leq \nu +L]}$ and $Z_{t} =(Z_{1t},...,Z_{Nt})^T$ are i.i.d. $N(0, \Sigma)$. $\nu$ is called the common change point with length $L$. 
By denoting $X_t = (X_{1t}, ..., X_{Nt})^T$, $\mu =(\mu_1,..., \mu_N)^T$, and $||\mu|| =(\mu^T\Sigma^{-1} \mu)^{1/2} $ as the signal strength, we can write $\mu =||\mu|| \Sigma^{-1/2}\gamma $ where $
\gamma = \Sigma^{-1/2} \mu /||\mu|| $ is the signal direction. Thus $X_t$ follows multivariate normal $N(||\mu||\Sigma^{1/2} \gamma  I_{[\nu \leq t \leq \nu +L]} , \Sigma )$. For observation $X_t$, given a reference value $||\delta||$ for the signal strength $||\mu||$, the log-likelihood ratio for $ \nu \leq t \leq \nu +L$ vs $t < \nu$ or $t>\nu +L$ is equal to
\begin{equation}
||\delta||\gamma^T \Sigma^{-1/2} X_t -\frac{||\delta||^2}{2} = ||\delta || (\gamma^T \Sigma^{-1/2} X_t -\frac{||\delta||}{2}) .
\end{equation}
In the following, we first give definitions for several chart that are direction invariant.

\noindent{\bf (1) Multivariate EWMA chart}

The first approach is to use EWMA of the log-likelihoods as discussed in Lowry, et al. (1992) and Wu and Wu (2022). 
By defining
\[
Y_{t}= (1-\beta)Y_{t-1} +\beta X_{t} 
\]
with $Y_0 =0$ for the same weight parameter $\beta$, the EWMA of the log-likelihoods (3.8) is equivalent to
\[
Y_{t}(\gamma) = \gamma^T\Sigma^{-1/2} Y_t .
\]
If we want the detection procedure is directional invariant, we use $
\max_{\gamma} Y_t(\gamma)$ as the detection process. This gives $\gamma=\Sigma^{-1/2}Y_t/(Y_t^T\Sigma^{-1}Y_t)^{1/2}$ and $
\max_{\gamma} Y_t(\gamma) =(Y_t^T\Sigma^{-1}Y_t)^{1/2}$. 

\noindent \underline{{\bf MEWMA procedure}}: {\em An alarm is made at
\begin{equation}
\tau_{MEW}  =\inf\{ t>0: ~  Y_{t}^T\Sigma^{-1}Y_t >b^2 (\frac{\beta}{2-\beta}) \},
\end{equation}
where the factor $\beta/(2-\beta)$ is the limiting variance for each component. }

\noindent{\bf (2) Multivariate MA chart}

The second approach is to take the moving average of the log-likelihoods with a window size $w$. By denoting $\bar{X}_{t;w} = \frac{1}{w} \sum_{t-w+1}^t X_j $, the moving average of the log-likelihoods can be written as 
\[
\gamma^T\Sigma^{-1/2}\bar{X}_{t;w} .
\]
It is obvious that $\gamma =\Sigma^{-1/2}\bar{X}_{t;w}/ (\bar{X}_{t;w}^T\Sigma^{-1} \bar{X}_{t;w})^{1/2} $ maximizes the above quantity. 

\noindent \underline{{\bf MMA procedure}:} {\em Make an alarm at 
\[
\tau_{MMA} =\inf\{t>w: (\bar{X}_{t;w}^T\Sigma^{-1} \bar{X}_{t;w})^{1/2}>h \}. 
\]
}

\noindent{\bf Remark :}  Just like the one-dimensional case, the $b^2 \beta /(2-\beta )$ plays the role for the target squared signal strength $||\delta||^2$ and $h$ plays the role for the target signal strength $||\delta||$. The design of $\beta$ and $w$ will also depend on the target length $L$.

\noindent{\bf (3) Multivariate CUSUM and GRLT charts}

There are many cases that a multivariate CUSUM procedure can be developed for detection a transient change as considered in Basseville and Nikiforov (1993, Section 7.2). Here we only consider the direction-invariant case. 

Note that the log-likelihood ratio of testing $H_0: ||\mu||=0$ vs $H_1: ||\mu|| =||\delta||$ for observations $X_{t-w+1},.., X_t$  can be written as
\begin{equation}
||\delta||\gamma^T  \Sigma^{-1/2}(X_{t-w+1}+...+X_t)-\frac{1}{2} w||\delta||^2  ,
\end{equation}
where $\gamma =\Sigma^{-1/2}\delta/ (\delta^T \Sigma^{-1} \delta)^{1/2}$. The $\gamma$ that maximizes the above log-likelihood ratio is equal to
\[
\Sigma^{-1/2}(X_{t-w+1}+...+X_t) / ((X_{t-w+1}+...+X_t)^T\Sigma^{-1}(X_{t-w+1}+...+X_t))^{1/2}.
\]
Thus, the maximum log-likelihood ratio is equal to
\begin{equation}
||\delta|| ((X_{t-w+1}+...+X_t)^T\Sigma^{-1}(X_{t-w+1}+...+X_t))^{1/2} -\frac{||\delta||^2}{2} w .
\end{equation}
By using the possible range of signal length as $[w_0, w_1]$, we can use the following windowed multivariate CUSUM procedure. 

\noindent \underline{{\bf MCUSUM Procedure:}} {\em An alarm is made at
\[
\tau_{MCU} =\inf\{t>w_1: \max_{w_0<w\leq w_1} w [(\bar{X}_{t;w}^T\Sigma^{-1} \bar{X}_{t;w})^{1/2} -||\delta||/2 ]>d\}
\]
}

If the signal strength $||\delta||$ is also unknown, we can estimate it from (3.11) as 
$ ((X_{t-w+1}+...+X_t)^T\Sigma^{-1}(X_{t-w+1}+...+X_t))^{1/2}/w $. So the following multivariate windowed-GLRT chart can be used:

\noindent \underline{{\bf Windowed GLRT Procedure:}} {\em Make an alarm at 
\[
\tau_{GLT} =\inf\{t>0: \max_{W_0\leq w <W_1} w (\bar{X}_{t;w}^T\Sigma^{-1} \bar{X}_{t;w}) >b^2\}
\]
}

\noindent{\bf Remark.} One can see that both MCUSUM and MGLRT charts are extensions of the MMA procedure. The design for the MCUSUM chart is more complicated as it also needs a reference value for the signal strength. If the signal length is $L$, then 
at $||\delta||$, $b \approx L ||\delta||/2$. For the MGLRT chart, if the target signal strength is $||\delta||$ with signal length $L$, then $b^2 \approx L ||\delta||^2$. 

\noindent{\bf Remark:} To mimic the CUSUM procedure in the one-dimensional case, one can update the common change point estimate at each time $t$ and use the cumulative sums after change point estimate in the recursive form (Pignatiello and Runger (1990)). 

{\em Let $\hat{\nu}_0 =0$ and for $t\geq 1$
\begin{equation}
MC1_t =max(0, ((\sum_{\hat{\nu}_t+1}^t X_i )^T \Sigma^{-1} (\sum_{\hat{\nu}_t+1}^t X_i ))^{1/2} - \frac{k_1}{2} (t- \hat{\nu}_t) )
\end{equation}
for a reference value $k_1 $ of the signal strength and reset $
\hat{\nu}_{t+1} =t, $
if $MC1_t =0$; otherwise $\hat{\nu}_{t+1}= \hat{\nu}_t$. An alarm will be made at
\[
\hat{\tau}_1 =\inf\{ t>0: MC1_t > h_1 \}.
\]
After detection, the change point can be estimated as $\hat{\nu}_{\hat{\tau}_1}$. }

\noindent{\bf Remark.} Monitoring charts based on window restricted  CUSUM and Shiryayev-Roberts rules by mixing likelihood ratio are studied in Siegmund and Yakir (2008) in terms of SADDT, where theoretical approximations for the FDP are also obtained.

\subsection{False Detection Probability} 

Without loss of generality, we can assume that $\Sigma= I_N$ is known. 
We  first give the approximation of FDP for MMA procedure.

\begin{theorem}
{\it Let $b= h \sqrt{w}$. As $w, N \rightarrow \infty$ such that $N/b^2 \rightarrow 0$. 
\[
P_0^*(\tau_{MMA} \leq L ) \approx \frac{2L}{w} \frac{(b^2/2)^{N/2}} {\Gamma(N/2) } e^{-b^2/2}  \upsilon (b\sqrt{2/w} ). 
\]
}
\end{theorem}

\noindent{\bf Proof:} First we can write 
\[
\tau_{MMA} =\inf\{t>w: \sum_{j=1}^N Z_{jt}^2 > h^2 w =b^2\},
\]
where $Z_{jt} = \frac{1}{\sqrt{w}} \sum_{t=t-w+1}^t X_{ji} $ which follows $N(0,1)$ under $P_0^*(.)$. Note that 
\[
\psi(\theta) = \ln E_0^*[\exp(\theta Z_{jt}^2 )] =\frac{1}{2}\ln (1-2\theta).
\]
Define a changed measure $P_t^*(.) $ on the time window $[t-w+1, t] $ such that
$X_{t-w-1}, ..., X_t$ are i.i.d. $N(0, 1/\sqrt{1-2\theta})$ for $\theta <1/2$ with changed variance instead of mean. This implies that $Z_{it}^2$ follows $\chi_1^2$ and the likelihood ratio for $X_1,,..., X_L$ has the form
\[
\frac{dP_t^*}{dP_0^*} (Z_{jt}^2) = \exp(\theta  Z_{jt}^2 +\frac{1}{2} \ln(1-2\theta )),
\]
for $j=1,2.,,.N$. 
Let $\theta $ satisfy $N \frac{d}{d\theta} [-\ln(1-2\theta)/2] =b^2$, i.e. $1- 2\theta =N/b^2$ or $\theta = (1-N/b^2)/2 $. By denoting
\[
l_t =\theta \sum_{j=1}^N Z_{jt}^2 +\frac{N}{2} \ln(1-2\theta );
\]
\[
\tilde{l}_t =\theta ( \sum_{j=1}^N Z_{jt}^2 -N/(1-2\theta)) ,
\]
and noting $\sigma_N^2 = Var_t^*(\tilde{l}_N ) = 2N\theta^2/(1-2\theta)^2$, we can approximate 
\[
P_0^*(\tau_{MMA} \leq L ) \approx 
L e^{-N(\theta/(1-2\theta) +\ln(1-2\theta)/2)} \frac{1}{\sqrt{2\pi}}\frac{1-2\theta}{\sqrt{2N} \theta }
\lim_{L\rightarrow \infty} E_t^*[ \frac{M_t}{S_t}]
\]
where $M_t =\max_{w<s\leq L} \exp (l_s-l_t) $ and $S_t = \sum _{w<s\leq L} \exp (l_s-l_t) $. 
By using Sterling's formula for large $z$
\[
\Gamma(z) =\sqrt{\frac{2\pi}{z}} (\frac{z}{e})^z (1+O(\frac{1}{z})),
\]
we can show that 
\[
 e^{-N(\theta/(1-2\theta) +\ln(1-2\theta)/2} \frac{1}{\sqrt{2\pi}}\frac{1-2\theta}{\sqrt{2N} \theta }
 \approx \frac{(b^2/2)^{N/2-1} e^{-b^2/2}}{\Gamma(N/2) (1-N/b^2)}.
 \]
 
 Finally, we evaluate the factor $E_t^*[ \frac{M_t}{S_t}]$. For $t-w \leq s \leq t$, we can write
 \begin{eqnarray*}
 Z_{js}^2-Z_{jt}^2 & = &(Z_{js}+Z_{jt})(Z_{js}-Z_{jt}) \\
 &= & \frac{1}{w} (\sum_{s-w+1}^{t-w}X_{ji} +2 \sum_{t-w+1}^{s}X_{ji} +\sum_{s+1}^t X_{ji} )(\sum_{s-w+1}^{t-w}X_{ji}-\sum_{s+1}^t X_{ji}) \\
 &=& \frac{t-s}{w}(Z_1 +\frac{Z_2}{\sqrt{1-2\theta}} +2 (\frac{w-(t-s)}{t-s})^{1/2} Z_3 ) (Z_1 -\frac{Z_2}{\sqrt{1-2\theta}} ) \\
 &= & \frac{t-s}{w} (Z_1^2 - \frac{Z_2^2}{1-2\theta} +2 (\frac{w-(t-s)}{t-s})^{1/2} Z_3  (Z_1 -\frac{Z_2}{\sqrt{1-2\theta}} )),
 \end{eqnarray*}
where $Z_1=\sum_{s-w+1}^{t-w}X_{ji}/\sqrt{t-s} $, $Z_2=\sum_{s+1}^t X_{ji}/\sqrt{t-s} $, and $Z_3=\sum_{t-w+1}^{s}X_{ji}/\sqrt{w-(t-s)} $ are independent standard normal random variables. 

Thus, 
\[
E_t^*(l_s -l_t ) = N\theta \frac{t-s}{w} (1-\frac{1}{1-2\theta}) = (t-s) \frac{-2\theta N}{w(1-2\theta )} =-\frac{b^2}{w} (1-\frac{N}{b^2}) (t-s).
\]
Also, as $w$ is large, by ignoring the lower order terms, we can get
\[
Var_t^*(l_s-l_t) \approx  N \theta^2 \frac{8(1-\theta)}{(1-2\theta)w} (t-s)
= \frac{2b^2}{w} (1-\frac{N}{b^2})^2 (1+\frac{N}{b^2}) (t-s).
\]
That means, $l_s -l_t$ can be approximated by a two-sided random walk with drift $-\frac{b^2}{w} (1-\frac{N}{b^2})$ and variance $\frac{2b^2}{w} (1-\frac{N}{b^2})^2 (1+\frac{N}{b^2})$. From Siegmund and Yakir (2000), we can approximate
\[
 E_t^*[ \frac{M_t}{S_t}] \approx \frac{b^2}{w} (1-\frac{N}{b^2}) \upsilon ( b(\frac{2}{w} (1-\frac{N}{b^2})^2 (1+\frac{N}{b^2}))^{1/2}).  \qed
 \]
 
Second, we give the approximation of FDP for MEWMA:

\begin{theorem}
As $b \rightarrow \infty$, $\beta \rightarrow 0$, and $N\rightarrow \infty$ such that $N/b^2 \rightarrow 0$ and $\beta b^2/N \rightarrow 0$, 
\[
P_0^*(\tau_{MEW}\leq L) \approx 2L\beta \frac{(b^2/2)^{N/2}}{\Gamma(N/2)} e^{-b^2/2} \upsilon (b\sqrt{2\beta}).
\]
\end{theorem}
\noindent{\bf Proof}. First, by denoting $Z_{jt} =Y_{jt}/\sqrt{\beta/(2-\beta)}$, we can write
\[
\tau_{MEW} =\inf\{ t>0: \sum_{j=1}^NZ_{jt}^2 > b^2 \} .
\]
Next, we define a changed measure $P_t^*(.)$ such that under $P_t^*(.)$, $X_0$ follows $N(0, \sqrt{\beta/(2-\beta)})$, $X_t $ follows $N(0,\frac{1}{\sqrt{\beta (2-\beta )}} (\frac{1}{1-2\theta} -(1-\beta)^2))^{1/2})$, and $X_i$ follow $N(0, 1) $ for $i \neq t$.   Then the likelihood ratio under $P_t^*(.)$ with respect to $P_0^*(.)$ is 
\[
\frac{dP_t^*}{dP_0*} (Z_{jt}^2) = \exp(\theta Z_{jt}^2 +\frac{1}{2}\ln (1-2\theta )).
\]
Select $\theta $ such that $N/(1-2\theta) =b^2$, So $1-2\theta = N/b^2$ or $\theta =(1-N/b^2)/2$.
Using the same notations as in the last theorem, we obtain
\[
P_0^*(\tau_{MEW} \leq t ) \approx L
\frac{(b^2/2)^{N/2-1} e^{-b^2/2}}{\Gamma(N/2) (1-N/b^2)} \lim_{L\rightarrow \infty} E_t^*[ \frac{M_t}{S_t}]
\]
where $M_t =\max_{0<s\leq L} \exp (l_s-l_t) $ and $S_t = \sum _{0<s\leq L} \exp (l_s-l_t) $. 

Finally, we evaluate $E_t^*[ \frac{M_t}{S_t}]$. We can write
\begin{eqnarray*}
Z_{j(t+1)}^2-Z_{jt}^2 &= &((1-\beta)Z_{jt}+\beta X_{j(t+1)}/\sqrt{\beta/(2-\beta)})^2 -Z_{jt}^2 \\
&\approx & -2\beta Z_{jt}^2 +2\sqrt{2\beta}Z_{jt} X_{j(t+1)}.
\end{eqnarray*}
Thus, the mean of $l_{t+1}-l_t$ is $2\theta N\beta /(1-2\theta) = (1-N/b^2) \beta b^2$ and variance is approximately equal to
\[
N\theta^2 (\frac{8\beta^2}{(1-2\theta)^2} +\frac{8\beta}{1-2\theta}) =2\beta b^2 (1-N/b^2)^2 (1+\beta b^2/N).
\]
Similarly, we can show the same results for $l_{t-1}-l_t$. That means
\[
E_t^*[ \frac{M_t}{S_t}] \approx (1-N/b^2) \beta b^2 \upsilon( b(1-N/b^2) \sqrt{2\beta (1+\beta b^2/N)}).
\qed
\]

The combination of the results in Theorems 1 and 5 gives the following approximation for windowed GRLT procedure:

\begin{cor} {\em As $b \rightarrow \infty$ and $L \rightarrow \infty$ such that $b/\sqrt{w_0}$ and $b/\sqrt{w_1}$ go to constants and $\ln(L)/b^2 \rightarrow 0$,
\[
P_0^*(\tau_{GLT} \leq L ) \approx \frac{(b^2/2)^{N/2-1} e^{-b^2/2}}{\Gamma(N/2)} \int_{w_0/L}^{w_1/L} \frac{b^4}{4t^2} \upsilon^2(\frac{b}{\sqrt{Lt}}) dt 
= 2L \frac{(b^2/2)^{N/2} e^{-b^2/2}}{\Gamma(N/2)} \int_{b/\sqrt{w_1}}^{b/\sqrt{w_0}} \frac{u \upsilon^2 (u)}{2} du. 
\]}
\end{cor}

Although we assume $N \rightarrow \infty$ above, when $N=1$, the results reduces to the one-dimensional case with two-sided signal. Finally, we give the result for the MCUSUM procedure and the proof is given in the Appendix 

\begin{theorem} {\it As $L, d \rightarrow \infty$ and $w_0, w_1 =O(d) $ such that $ \frac{d}{w_1} < \frac{||\delta||}{2} < \frac{d}{w_0}$ and $Le^{-||\delta|| d} \rightarrow 0$,
\begin{eqnarray*}
P_0^*(\tau_{MCU} \leq L) &\approx & \frac{L||\delta||^2}{2} (4||\delta||d)^{(N-1)/2} \frac{\Gamma((N-1)/2)}{\Gamma(N-1)} e^{-||\delta||d} \upsilon^2(||\delta||) \\
&\approx & \frac{L||\delta||^2}{2} (4||\delta||d)^{(N-1)/2} \frac{\Gamma((N-1)/2)}{\Gamma(N-1)}e^{-||\delta||(d +2\rho_+) }. 
\end{eqnarray*}
}
\end{theorem}

In Table 2, we report some simulation comparisons of simulated FDP and the approximations (given in brackets) given in the above theorems for $N=20$. For the approximated values, we use the simple approximation $\upsilon (x) =exp(-0.5826 x)$ for the numerical computation. 
 We see the approximations of FDP for MMA and MEWMA are very satisfactory. For the GLRT, the approximation gets better when the tail area is smaller, as expected. 

\begin{table}[ht]
\begin{center}
Table 2. Comparison of approximated and simulated FDP \\
\begin{tabular}{|c |c |c | c|c |}   \hline 
b & EWMA & MA & MA &  GLRT\\
 & $\beta=0.05$ & $w=10$ &$w=20$ & $20<w\leq 50$ \\ \hline 
6.0 & 0.0985(0.0992) & 0.1466(0.1324) & 0.0983(0.0998)  & 0.2456(0.4076) \\
6.1 & 0.0724(0.7403) &0.1106(0.0984)& 0.0770(0.0746) &  0.1916(0.3097) \\
6.2 & 0.0543(0.5440) & 0.0846(0.0720) & 0.0565 (0.0549)  & 0.1483(0.2318) \\
6.3 & 0.0400(0.0394)& 0.0598(0.0519) & 0.0393 (0.0398)  & 0.1162 (0.1709) \\
6.4 & 0.0274(0.0281)& 0.0422(0.0368) & 0.0288 (0.0284)  & 0.0842(0.1240) \\
6.5 & 0.0190(0.0197) & 0.0308(0.0258) & 0.0205(0.0200) & 0.0556(0.0887) \\
6.6 & 0.0140 (0.0136)& 0.0207(0.0178) & 0.0143(0.0138)  & 0.0436 (0.0625)\\
6.7 & 0.0096(0.0093)   & 0.0145(0.0121) &0.0097(0.0095)  & 0.0284 (0.0434)\\
6.8 & 0.0065(0.0063) & 0.0093(0.0081)& 0.0074(0.0066)&0.0226(0.0298) \\ 
6.9 &  0.0044(0.0041)   & 0.0066(0.0054)  & 0.0043(0.0042) & 0.0178 (0.0200) \\
7.0 &  0.0026(0.0027)   & 0.0041(0.0035) & 0.0030(0.0028)  &  0.0102(0.0134)\\  \hline
\end{tabular}
\end{center}
\end{table}

\subsection{Comparison of POD}

For a complete comparison of POD, we consider two extreme scenarios. First, we assume that the signal appears at equal strength in all channels. 

In Table 3, for N=20, we take $FDP =0.02$ and signal $||\mu|| (1, ..., 1 )^T$ for $||\mu|| =$0.1, 0.2, 0.25, 0.3, 0.4, 0.5 with strength $||\mu|| \sqrt{20}$.\\
(i) For MEWMA chart, $b =6.5$ and $\beta =0.05$; \\
(ii) For MMA chart , we take $w=10, 20, 50 $ that correspond to $b=6.6, 6.5, 6.37$ respectively;\\
(iii) For the windowed GLRT chart, we take $b=6.84$ with $20<w\leq 50$;\\
(iv) For the windowed CUSUM chart, we take $b=20.3 $ with reference value for signal strength  $0.25 \sqrt{20}$ and $20<w\leq 50$. 

The simulations are replicated 50,000 times for EWMA and MA charts and 5,000 for MCUSUM and GLRT charts. The comparison shows that MA procedure performs better when the signal length matches the window size, as expected. However, when signal length is longer than the window width,  the MA performs poor compared with other procedures. The GLRT procedure can compromise this shortcoming when the signal length is within the window. On the other hand, when the signal length gets longer, the EWMA procedure performs better than other procedures. 

\begin{table}[ht]
\begin{center}
Table 3. Comparison of POD for $N=20$ with all channel change\\
\begin{tabular}{|c|c |c |c |c|  c| c|c |}   \hline
L &$||\mu||$ & EWMA & MA & MA & MA & CUSUM & GLRT\\
& parameter & $\beta=0.05$ &$w=10$ &$w=20$ & $w=50$& $20<w\leq 50$  & $20<w\leq 50$ \\
 &boundary b & 6.5 & 6.6 & 6.5 & 6.37 & d=20.3   & 6.84  \\ \hline
20 &0.0 & 0.0198 & 0.0219& 0.0209 &0.0199& 0.0228 & 0.0195 \\
  &0.1 & 0.0441 &0.0446 & 0.0507 &0.0276 & 0.0316 & 0.0386 \\
  &0.2& 0.2547 & 0.2229 & 0.3250 &0.0806& 0.2104 & 0.2204\\
  &0.25& 0.5037 & 0.4367 & 0.6280 &0.1463&  0.4536 & 0.5024 \\
  &0.3 & 0.7667 & 0.6921 & 0.8743 &0.2763&  0.7378 & 0.7754 \\
  &0.4 & 0.9875 & 0.9760& 0.9983 &0.6572&  0.9894 & 0.9920 \\
  &0.5 & 0.99996 & 0.99984 & 1.00 & 1.00 & 1.00 & 1.00\\ \hline
30 &  0.0 & 0.0287 & 0.0323 & 0.0301 &0.0283 & 0.0272 & 0.0256 \\
  &0.1 & 0.0920 &0.0706 & 0.0926 &0.0547& 0.0864& 0.0842 \\
  &0.2& 0.5643 & 0.3449 & 0.5541&0.2973&  0.6014 & 0.5906\\
  &0.25& 0.8582 & 0.6180 & 0.8467 &0.5656 & 0.8842 & 0.8798 \\
  &0.3 & 0.9795 & 0.8635 & 0.9752 &0.8381& 0.9894 & 0.9868\\
  &0.4 & 0.99994 & 0.99768 & 0.99994 &0.9959& 1.0000 & 1.0000 \\
  &0.5 & 1.0 & 1.00& 1.00 &1.00 & 1.00 & 1.00\\ \hline
50 &  0.0 & 0.0460 & 0.0513 & 0.0485 &0.0438& 0.0396 & 0.0424 \\
  &0.1 & 0.2154 &0.1232 & 0.1740 &0.2133& 0.2060 & 0.1970 \\
  &0.2& 0.9049 & 0.5344 & 0.8035 &0.9517& 0.9284 & 0.9214\\
  &0.25& 0.9948 & 0.8264 & 0.9733 &0.9989& 0.9978 & 0.9972 \\
  &0.3 & 0.99992 & 0.97066 & 0.99914 &1.0000 &1.0000 & 1.0000\\
  &0.4 & 1.00 & 1.0 & 1.00 &1.00 & 1.00 & 1.00 \\
  &0.5 & 1.00 & 1.00& 1.00 & 1.00 &1.00  & 1.00 \\ \hline
\end{tabular}
\end{center}
\end{table}

Second, we consider another extreme case that the signal only appears in one channel, i.e $\mu =||\mu|| (1, 0,..., 0)^T$ with $||\mu|| =0.0, 0.25, 0.5, 0.75, 1.0,1.25, 1.5, 1.75, 2.0$.  Table 4 gives the corresponding results. 

\begin{table}[ht]
\begin{center}
Table 4.  Comparison of POD for $N=20$ with only one channel change \\
\begin{tabular}{|c|c |c |c | c|c |c|c |}   \hline
L & $||\mu||$ & EWMA & MA & MA & MA &CUSUM & GLRT\\
 &parameter & $\beta=0.05$ & $w=10$ &$w=20$ & $w=50$& $20<w\leq 50$  & $20<w\leq 50$ \\
 & boundary b & 6.5 & 6.6 & 6.5 & 6.37 & d=20.3   & 6.84  \\ \hline
20 &0.0 & 0.0195 & 0.0206& 0.0218 &0.194&  0.0202 & 0.0206 \\
 &0.25 & 0.0252 &0.0268& 0.0270 & 0.0204& 0.0222  & 0.0218\\
 &0.5 & 0.0539 & 0.0541 & 0.0622 & 0.0331& 0.0436  & 0.0408\\
 & 0.75& 0.1432 & 0.1326 & 0.1809 &0.0525 &  0.1146 & 0.1208 \\
 & 1.0 & 0.3582 & 0.3116 & 0.4603 & 0.1066 & 0.3022 & 0.3370\\
 &1.25 & 0.6642 & 0.5854 & 0.7928 &0.2211&  0.6430 & 0.6690\\
 &1.5 & 0.8973 & 0.8434 & 0.9599 & 0.3931 & 0.8884 & 0.9076\\ 
 &1.75& 0.9835&  0.9678     &  0.9968     & 0.6252 &  0.9874   & 0.9878    \\
 &2.00& 0.99872  & 0.9972 &  .99988      & 0.8086 &  0.9986    & 0.9994    \\\hline
30 &0.0 & 0.0271 & 0.0315 & 0.0305 &0.0264& 0.0274 & 0.0288 \\
 &0.25 & 0.0426 &0.0404& 0.0455 &0.0349 & 0.0406 & 0.0428 \\
 &0.5 & 0.1177 & 0.0840 & 0.1196 &0.0606& 0.1112  & 0.1060\\
 &0.75& 0.3539 & 0.2119 & 0.3507 &0.1831 & 0.3566 & 0.3344 \\
 &1.0 & 0.7212 & 0.4675 & 0.7111 &0.4276 & 0.7468 & 0.7438 \\
 &1.25 & 0.9490 & 0.7753 & 0.9438 &0.7400 & 0.9674 & 0.9622 \\
 &1.5 & 0.9965& 0.9530 & 0.9960 &0.9394 & 0.9982 & 0.9980\\ 
 &1.75& 0.99988 &  0.9965    &  0.99992  &0.9926 &  1.00   & 1.00    \\
 & 2.00& 1.00   & 0.99992      & 1.00     &0.9998 &  1.00&  1.00   \\\hline 
50 & 0.0 & 0.0480 & 0.0513 & 0.0484 &0.0421&  0.0398 & 0.0432 \\
 &0.25 & 0.0799 &0.0682 & 0.0756 &0.0711 & 0.0732 & 0.0672 \\
 &0.5 & 0.2814 & 0.1445 & 0.2210 &0.2966 & 0.2702  & 0.2642\\
 &0.75& 0.7116 & 0.3471 & 0.5727 &0.7862 & 0.7574 & 0.7316 \\
 &1.0 & 0.9693 & 0.6835 & 0.9117 &0.9876 & 0.9798 & 0.9802 \\
 &1.25 & 0.99952  & 0.9330 & 0.9959 &0.9999& 1.0000 & 0.9994 \\
 &1.5 & 1.00 & 0.9963 & 1.00 &1.00& 1.00  & 1.00\\ 
 &1.75& 1.00   & 0.99998  & 1.00 &1.00     &  1.00    &  1.00    \\
 & 2.00& 1.00  & 1.0     & 1.00 &1.00     &  1.00     & 1.00    \\\hline
\end{tabular}
\end{center}
\end{table}

Comparing Tables 4 with Tables 1 we see that if the signal only appears in one channel, the power of detection is reduced dramatically for the multivariate monitoring chart. In next section, we shall consider some techniques on how to improve the POD.

In the following, we present some theoretical results of POD. 

(1) For the MEWMA chart, We shall make the assumption that $\beta \rightarrow 0$, $b \rightarrow \infty$, and $L \rightarrow \infty$ such that $ b\sqrt{\beta/(2-\beta)} \rightarrow h$, and $L\beta \rightarrow \infty$. Under $P_{\delta}^*(.)$, we can write 
\[
P_{\delta}^*(\tau_{MEW} \leq L ) 
=P_0^* (\max_{0\leq t \leq L} (Y_t^TY_t +2(1-e^{-\beta t}) \delta^T Y_t +(1-e^{-\beta t})^2 \delta^T\delta)>h^2).
\]
By ignoring the overshoot at the boundary crossing time, we have
\[
Y_t^TY_t +2(1-e^{-\beta t}) ||\delta|| \gamma^T Y_t +(1-e^{-\beta t})^2 \delta^T\delta =h^2,
\]
where $\gamma =\delta /||\delta||$. Similar to the one-dimensional case, when $||\delta|| <h$, since the middle term on the left side has mean zero, we can ignore it. By denoting 
\[
b(t) =(b^2 -(1-e^{-\beta t})^2 ||\delta||^2 (2-\beta)/\beta)^{1/2}, 
\]
we have the following approximation for the local POD:
\[
P_{\mu}^*(\tau_{MEW}\leq L) \approx 2 \int_0^{L\beta} \frac{(b^2(t)/2)^{N/2}}{\Gamma(N/2)} e^{-b^2(t)/2} 
\upsilon (b(t) \sqrt{2\beta} ) dt. 
\]

When $||\delta|| \geq h$, we have
\begin{eqnarray*}
1-e^{-\beta t} &=&\frac{h}{||\delta||} (1-\frac{2||\delta||}{h^2} (1-e^{-\beta t}) \gamma^TY_t -\frac{1}{h^2} Y_t^TY_t )^{1/2} \\
&= &\frac{h}{||\delta||} (1-\frac{||\delta||}{h^2}(1-e^{-\beta t}) \gamma^TY_t 
-\frac{1}{2h^2} Y_t^TY_t -\frac{||\delta||^2}{2h^4}(1-e^{-\beta t})^2 
(\gamma^TY_t)^2 +o_p(\beta)).
\end{eqnarray*}
From above, we have
\[
1-e^{-\beta t}=\frac{h}{||\delta||} - \frac{1}{h} (1-e^{-\beta t}) \gamma^TY_t +o_p(\sqrt{\beta}) 
\]
\[
=\frac{h}{||\delta||} -\frac{1}{||\delta||} \gamma^TY_t +o_p (\sqrt{\beta}).
\]
By re-plugging back, we have
\[
1-e^{-\beta t} =\frac{h}{||\delta||}-\frac{1}{||\delta||} \gamma^TY_t +
\frac{1}{2h||\delta||} (\gamma^TY_t)^2 -\frac{1}{2h||\delta|| } Y_t^TY_t + o_p(\beta) . 
\]
From this, we get
\begin{eqnarray*}
-t\beta &=& \ln (1-\frac{h}{||\delta||} +\frac{1}{||\delta||} \gamma^TY_t
-\frac{1}{2h||\delta||} (\gamma^TY_t)^2 +\frac{1}{2h||\delta||} Y_t^TY_t +o_p(\beta)) \\
&=& \ln (1-\frac{h}{||\delta||}) + \frac{1}{||\delta||-h} \gamma^TY_t -
\frac{||\delta||}{2h(||\delta||-h)^2} (\gamma^TY_t)^2 + \frac{1}{2h(||\delta||-h)} Y_t^TY_t +o_p(\beta) .
\end{eqnarray*}
Therefore, we have
\begin{eqnarray*}
E_{\delta}^*(\tau_{MEW}) &=& -\frac{1}{\beta} \ln (1-\frac{h}{||\delta||}) +
\frac{||\delta||}{4h(||\delta||-h)^2}
-\frac{N}{4h(||\delta||-h)} +o(\beta), \\
Var_{\delta}^*(\tau_{MEW}) &=& \frac{1}{2\beta (||\delta||-h)^2} +O(1)
\end{eqnarray*}
An approximation for POD can be obtained by the normal law as
\begin{equation}
P_{\delta}^*(\tau_{MEW} \leq L ) \approx \Phi\left ( \frac{\beta L+\ln (1-h/||\delta||) -\beta ||\delta||/(4h(||\delta||-h)^2) +N\beta/(4h(||\delta||-h))}{\sqrt{\beta/2}/(||\delta||-h) } \right )
\end{equation}

(2) For the MMA chart, we assume $w \rightarrow \infty$ and $ L/w \rightarrow \infty$. Note that under $P_{\delta}^*(.)$ at the boundary crossing $t$, by ignoring the overshoot, we have
\[
(\bar{X}_{t;w} +\delta \min(t/w,1))^T(\bar{X}_{t_w} +\delta \min (t/w,1)) =\bar{X}_{t;w}^T \bar{X}_{t;w} +2 (\min (t/w,1)) \delta^T\bar{X}_{t;w} +(\min(t/w,1))^2 \delta^T\delta^2 = h^2.
\]
For $||\delta|| < h$, we can ignore the middle term, and define
\[
h(t)=(h^2-\delta^T\delta \min(t^2 ,1))^{1/2}.
\]
The local power can be approximated as
\[
P_{\delta}^*(\tau_{MMA} \leq L) \approx 2 \int_0^{L/w} \frac{(h^2(t) w/2)^{N/2}}{\Gamma(N/2)} e^{-h^2(t) w/2} \upsilon( h(t) \sqrt{2} ) dt.
\]
For $||\delta|| \geq h$, we have 
\begin{eqnarray*}
\frac{t}{w}& = &  \frac{h}{||\delta||} (1-\frac{2t}{h^2w} \delta^T\bar{X}_{t;w} -\frac{1}{h^2} \bar{X}_{t;w}^T \bar{X}_{t;w})^{1/2} \\
&= & \frac{h}{||\delta||}(1-\frac{t}{h^2w} \delta^T\bar{X}_{t;w} -\frac{1}{2h^2} \bar{X}_{t;w}^T \bar{X}_{t;w}
-\frac{(t/w)^2}{2h^4}( \delta^T\bar{X}_{t;w}
)^2 +o_p(1/w)) \\
&=& \frac{h}{||\delta||}- \frac{t/w}{h||\delta||}\delta^T\bar{X}_{t;w} -\frac{1}{2h||\delta||} \bar{X}_{t;w}^T \bar{X}_{t;w}-\frac{(t/w)^2}{2h^3||\delta||}( \delta^T\bar{X}_{t;w}
)^2 +o_p(1/w).
\end{eqnarray*}
By plugging in 
\[
\frac{t}{w}=\frac{h}{||\delta||}- \frac{t/w}{h||\delta||}\delta^T\bar{X}_{t;w} +o_p(1/\sqrt{w})
\]
back in the above expression, we have
\begin{eqnarray*}
\frac{t}{w}&=&\frac{h}{||\delta||}- \frac{1}{h||\delta||} ( \frac{h}{||\delta||}- \frac{t/w}{h||\delta||}\delta^T\bar{X}_{t;w})\delta^T\bar{X}_{t;w} -\frac{1}{2h||\delta||} \bar{X}_{t;w}^T \bar{X}_{t;w}-\frac{1}{2h||\delta||^3}( \delta^T\bar{X}_{t;w}
)^2 +o_p(1/w) \\
&=&\frac{h}{||\delta||}- \frac{1}{||\delta||^2}\delta^T\bar{X}_{t;w} +\frac{1}{2h||\delta||^3}( \delta^T\bar{X}_{t;w}
)^2 -\frac{1}{2h||\delta||} \bar{X}_{t;w}^T \bar{X}_{t;w} +o_p(1/w) .
\end{eqnarray*}
This gives
\begin{eqnarray*}
E_{\delta}^*(\tau_{MMA}) &=& \frac{hw}{||\delta||} -\frac{N-1}{2h||\delta|| } +o(1), \\
Var_{\delta}^*(\tau_{MMA}) &= &\frac{w}{||\delta||^2 } +O(1).
\end{eqnarray*}
Using the normal law, the POD can be approximated as
\begin{equation}
P_{\delta}^*(\tau_{MMA}\leq L) \approx \Phi \left(\sqrt{w}||\delta|| \left (\frac{L}{w} -\frac{h}{||\delta||} +\frac{N-1}{2h||\delta|| w}\right)\right).
\end{equation}

(3) The approximation of POD for the MCUSUM chart can be obtained similar to the one-dimensional case. Assume $||\mu >||\delta||/2$ and $d/w_1 < ||\mu||-||\delta||/2 <d/w_0 $ . At the alarm time, we can write
\begin{eqnarray*}
& &  \max w(((\frac{Z_{t;w}}{\sqrt{w}} +\mu \min (\frac{t-\nu}{w},1))^T(\frac{Z_{t;w}}{\sqrt{w}} +\mu \min (\frac{t-\nu}{w},1)))^{1/2} -||\delta||/2 )  \\
& &  =\max((\mu^T\mu (\min(\frac{t-\nu}{w},1))^2 +2 \min(\frac{t-\nu}{w},1) \frac{\mu^TZ_{t;w}}{\sqrt{w}} +\frac{Z_{t;w}^TZ_{t;w}}{w})^{1/2} -||\delta||/2) =d .
\end{eqnarray*}
As $d \rightarrow \infty$, $w \rightarrow \infty$ and the maximum $w$ at the first order is $w=t-\nu$. A second order expansion gives
\[
(\mu^T\mu +2 \frac{\mu^TZ_{t;w}}{\sqrt{w}} +\frac{Z_{t;w}^TZ_{t;w}}{w})^{1/2}
=||\mu|| + \frac{\mu^TZ_{t;w}}{\sqrt{w} ||\mu||} + \frac{1}{2||\mu|| w} (Z_{t;w}^TZ_{t;w} - (\mu^TZ_{t;w})^2/||\mu||^2 ) +o_p(1/w). 
\]
So the above equation reduces to
\begin{eqnarray*}
d & =& (t-\nu) (||\mu|| -||\delta||/2) +\frac{\mu^TZ_{t;w}}{||\mu||} \sqrt{t-\nu} +
\frac{1}{2||\mu||} (Z_{t;w}^TZ_{t;w} - (\mu^TZ_{t;w})^2/||\mu||^2 ) +o_p(1) \\
&=& (t-\nu) (||\mu|| -||\delta||/2) +\frac{\mu^TZ_{t;w}}{||\mu||} \frac{\sqrt{d}}{(||\mu||-||\delta||/2)^{1/2}} (1- \frac{\mu^TZ_{t;w}}{\sqrt{d} ||\mu|| \sqrt{||\mu||-||\delta||/2}} )^{1/2} \\
& & +
\frac{1}{2||\mu||} (Z_{t;w}^TZ_{t;w} - (\mu^TZ_{t;w})^2/||\mu||^2 ) +o_p(1) \\
&=& (t-\nu) (||\mu|| -||\delta||/2) +\frac{\mu^TZ_{t;w}}{||\mu||} \frac{\sqrt{d}}{(||\mu||-||\delta||/2)^{1/2}} -\frac{(\mu^TZ_{t;w})^2}{2||\mu||^2 (||\mu|| -||\delta||/2)}\\
& & +
\frac{1}{2||\mu||} (Z_{t;w}^TZ_{t;w} - (\mu^TZ_{t;w})^2/||\mu||^2 ) +o_p(1).
\end{eqnarray*}
From this, we get
\begin{eqnarray*}
E_{\mu}^*(\tau_{MCU}) &= &\frac{d}{||\mu|| -||\delta||/2}+\frac{1}{2(||\mu|| -||\delta||/2)^2} -\frac{N-1}{2||\mu|| (||\mu|| -||\delta||/2)} +o(1), \\
Var_{\mu}^*(\tau_{MCU}) &= & \frac{d}{(||\mu|| -||\delta||/2)^3} +O(1). 
\end{eqnarray*}
A norm approximation can be used after continuous correction as
\[
P_{\mu}^*(\tau_{MCU} \leq L) =\Phi \left(\frac{L-d/(||\mu|| -||\delta||/2) -1/2(||\mu|| -||\delta||/2)^2 +(N-1)/2(||\mu|| -||\delta||/2)}{\sqrt{d}/(||\mu|| -||\delta||/2)^{3/2}} \right).
\]

(4) The derivation for the GLRT chart is straight forward. Assume $ b/\sqrt{w_1} < ||\mu|| <b/\sqrt{w_0} $. Again, at the alarm time with signal $\mu$, we have
\[
\max w(||\mu||^2 \min((t-\nu)^2/w^2,1) +2\min((t-\nu)/w,1) \frac{\mu^TZ_{t;w}}{\sqrt{w}} +\frac{Z_{t;w}^TZ_{t;w}}{w} ) = b^2.
\]
As $b\rightarrow \infty$, at $w=t-\nu$, we have
\[
||\mu||^2 (t-\nu) +2\sqrt{t-\nu} \mu^TZ_{t;w} +Z_{t;w}^TZ_{t;w} =b^2. 
\]
The same technique by plugging in the expansion for $\sqrt{t-\nu}$ gives
\begin{eqnarray*}
t-\nu &=& \frac{1}{||\mu||^2} (b^2 -2\sqrt{t-\nu} \mu^TZ_{t;w} -Z_{t;w}^TZ_{t;w} ) \\
&= &\frac{b^2}{||\mu||^2} -\frac{2}{||\mu||^2} (\frac{b}{||\mu||} - \mu^T Z_{t;w} ) \mu^TZ_{t;w} - \frac{Z_{t;w}^TZ_{t;w}}{||\mu||^2} +o_p(1) \\
&= &\frac{b^2}{||\mu||^2} -\frac{ 2b}{||\mu||^3} \mu^TZ_{t;w} +\frac{2(\mu^TZ_{t;w})^2}{||\mu||^2} +o_p(1).
\end{eqnarray*}
From above, we get
\[
E_{\mu}^*(\tau_{GLR}) = \frac{b^2+N}{||\mu||^2} +2 +o(1), 
Var_{\mu}^*(\tau_{GLR}) =  \frac{4b^2}{||\mu||^4} +O(1).
\]
So the normal approximation is 
\[
P_{\mu}^*(\tau_{GLR} \leq L) \approx \Phi \left (\frac{L- (b^2+N)/||\mu||^2 -2}{2b/||\mu||^2} \right).
\]

\section{Improving POD for a Portion of Changed Channels and a Real Example}

All the multivariate charts discussed above are proposed as change direction invariant. In signal processing, one may be more interested in the situation when only a small portion of channels have changes with similar strength, Siegmund, et al. (2011) proposed two soft-threshold approaches: one is the mixed likelihood ratio statistic and the other is the weighted sum of chi-square statistics, both can be treated as special functions of chi-square statistics that rely on a reference value for the changed proportion $p$. The third one is the hard-threshold by truncating off the smaller chi-square values that is larger than a certain reference value for the strength of signal. 

Following Table 3, we compare MEWMA, MMA, and GLRT charts by using threshold method. The first is the weighted MEWMA chart by using the weights proposed by Siegmund et al. (2011) and an alarm is made at
\begin{equation}
\tau_1 =\inf\{ t>0: \sum_{j=1}^N \frac{e^{Y_{jt}^2/2}}{(1-p)/p+ e^{Y_{jt}^2/2}} Y_{jt}^2> d^2 \}. 
\end{equation}

Here we take $p=0.1$ so $(1-p)/p=9$. The second is the hard-threshold MEWMA chart by dropping the channels with $|Y_{it}|^2 \leq 0.25$. So we make an alarm at
\begin{equation}
\tau_2 =\inf\{t>0: \sum_{j=1}^N Y_{jt}^2I_{[|Y_{jt}|>0.5]} > b^2 \frac{\beta}{2-\beta}. \}
\end{equation}
The third and the fourth are the MMA charts with window 20 and 10 and hard-threshold 0.25 for the squared average. The boundary is denoted as $h^2$. 
The last is the windowed GLRT with the same hard-threshold for the mean square of each channel.  

The simulation results from Table 5 show that for the EWMA procedure, the hard-threshold performs not only better than the weighted chart, but also better than the MA and GLRT charts for relatively large signal strength. Overall, the threshold method can improve the POD by 20\% to 50\% by comparing Table 5 with Table 4. The simulation results also show that unless the changed proportion is small, say  10\% to 15\% , the improvement by using the above techniques are insignificant. 

An alternative approach is to use the average of chi-square values that is larger than the threshold. The simulation results show similar results. 

\begin{table}[ht]
\begin{center}
Table 5.  Comparison of POD for $N=20$ with only one channel change \\
\begin{tabular}{|c|c |c |c | c| c|c |}   \hline
L &$||\mu||$ & EWMA(weighted) & EWMA(0.25) & MA(0.25) & MA(0.25) & GLRT(0.25)\\
 &parameter & $\beta=0.05$ & $\beta=0.05$ &$w=20$ & $w=10$  & $20 <w \leq 50$ \\
 &boundary  & $\frac{b^2\beta}{(2-\beta)}=$0.1165 & 0.396 & $h^2=1.26$ & 3.47  & $b=5.12$  \\ \hline
10 & 0.0 & 0.0113 & 0.0106& 0.0111 & 0.0107 & 0.0126\\
 &0.25 & 0.0120 &0.0117& 0.0121 & 0.0126  & 0.0132\\
 &0.5 & 0.0162 & 0.0180 & 0.0155 & 0.0199 & 0.0174\\
 &0.75& 0.0279 & 0.0406 & 0.0233 & 0.0413 & 0.0260 \\
 &1.0 & 0.0658 & 0.1081 & 0.0418 & 0.1003 & 0.0400 \\
 & 1.25 & 0.1442 & 0.2503 & 0.0829 & 0.2358 & 0.0728 \\
 & 1.5 & 0.2839 & 0.4685 & 0.1600 & 0.4519 & 0.1500\\ 
 & 1.75& 0.4858 &  0.7003 &  0.2829&  0.6978 & 0.2630 \\
 & 2.00& 0.6932 & 0.8682 & 0.4514& 0.8843 & 0.4210 \\\hline
20 &0.0 & 0.0191 & 0.0190& 0.0213 & 0.0192 & 0.0234\\
 &0.25 & 0.02466 &0.0241& 0.0272 & 0.0311  & 0.0284\\
 &0.5 & 0.05486 & 0.0695 & 0.0671 & 0.0516 & 0.0588\\
 &0.75& 0.1731 & 0.2764 & 0.2067 & 0.1262 & 0.1922 \\
 &1.0 & 0.4338 & 0.6217 & 0.5167 & 0.3000 & 0.4722 \\
 &1.25 & 0.7607 & 0.9028 & 0.8433 & 0.5623 & 0.7946 \\
 &1.5 & 0.9438 & 0.9870 & 0.9767 & 0.8265 & 0.9680\\ 
 &1.75& 0.9940 &  0.9993   &  0.9987     &  0.9576    & 0.9972  \\
 &2.00& 0.99968  & 1.00 & 0.9999      &   0.9960  & 1.00   \\\hline 
30 &0.0 & 0.0276 & 0.0265 & 0.0289 & 0.0343 & 0.0314 \\
 & 0.25 & 0.0401 &0.0414& 0.0429 & 0.0411 & 0.0488 \\
 &0.5 & 0.1292 & 0.1783 & 0.1270 & 0.0807  & 0.1524\\
 &0.75& 0.4129 & 0.5957 & 0.3928 & 0.1980 & 0.4520 \\
 &1.0 & 0.8015 & 0.9248 & 0.7569 & 0.4535 & 0.8474 \\
 &1.25 & 0.9744 & 0.9965 & 0.9632 & 0.7512 & 0.9898 \\
 &1.5 & 0.9993 & 1.00 & 0.9984 & 0.9456 & 0.9998\\ 
 &1.75& 0.99998 &  1.00   &  1.00    &  0.9960    & 1.00   \\
 & 2.00& 1.00   &  1.00     & 1.00     &0.9998 &  1.00   \\\hline 
50 &0.0 & 0.0447 & 0.0464 & 0.0474 & 0.0512 & 0.0472 \\
 &0.25 & 0.0825 &0.0929 & 0.0772 & 0.0673  & 0.0790 \\
 &0.5 & 0.3168 & 0.4436 & 0.2342 & 0.1424  & 0.3148\\
 &0.75& 0.7840 & 0.9094 & 0.6427 & 0.3355 & 0.8030 \\
 &1.0 & 0.9869 & 0.9977 & 0.9400 & 0.6590 & 0.9944 \\
 &1.25 & 0.9998  & 1.00 & 0.9983 & 0.9238 & 1.00 \\
 & 1.5 & 1.00 & 1.00 & 1.00  &0.9949  & 1.00  \\ 
 &1.75& 1.00   & 1.00  & 1.00      &  1.00    &  1.00    \\
 & 2.00& 1.00  & 1.00    & 1.00      &  1.00     & 1.00    \\\hline
\end{tabular}
\end{center}
\end{table}

\noindent{\bf Example.} We use Dow Jones 30 industrial stock prices (downloaded from finance.yahoo.com) as illustration of using the EWMA chart. The simple geometric random walk model is applied as the ACF of the differences of logarithm and its square do not show significant correlations.    

For the single sequence, the first figure of Figure 1 gives the plot of closing price of CVX from May 6, 2021 to May 6, 2022. The second figure gives the differences of logarithm. The third figure gives the plot after standardization by dividing its standard deviation. The fourth figure gives the plot of EWMA chart with $\beta =0.05$. By using the same design as given in Table 1 with FDP=0.01 for $L=20$ with limit 0.473, a transient signal is detected at $t=208$, continuing until 215. 

For multi-dimensional data stream, to match with the results given in Tables 3 to 5, we choose 20 of the 30 Dow Jones stocks. 
Figure 2 gives the plots of three MEWMA charts with $\beta =0.05$ for the first 20 stocks of DOW JONES 30:\\
MMM, AXP, APPL, BA, CAT, CVX, CSCO, DOW, GS, HON, \\
IBM, INTC, JNJ, JPM, MCD, MRK, MSFT, NKE, PG, CRM.

All stock prices are standardized as for CVX. The first plot is the MEWMA chart by using the estimated correlation matrix for the 20 stocks using the whole year's data as defined in (3.9). We found that the stocks in the same sectors have correlations larger than 0.5. For example, AXP, BA, GS, and JPM in the Finance Sector have high correlations, DOW, HON, and MMM are highly correlated, and APPL, MSFT, INTC, CSCO, and CRM are highly correlated. Otherwise, most correlations can be ignored. 
The largest eigenvalue is calculated as 7.7 and others are insignificant. This suggests that we can use block intra-class matrix to model the correlation matrix. 

The second figure gives the plot for the MEWMA chart by ignoring the correlation, i.e. treating observations from different channels as independent. The third chart gives the plot of MEWMA by using the hard threshold 0.5 from the second chart. 

We find that the first two charts behaves similarly and raise alarm at day 201 and 202, respectively (with limit 1.0833 from Table 4). The first chart shows an obvious outlier at day 145 that may be caused by the strong correlation from several stocks in the same sector. The third figure with hard threshold also raises an alarm at day 210 (with limit 0.396 from Table 5) similar to the second figure. By checking all 20 individual EWMA charts, we find that only CVX shows significant mean increase that causes alarm for the multivariate EWMA chart.  

In summary, we find that it is convenient to monitor both the individual and global EWMA charts simultaneously to detect the common signal that may appear in one or more channels and also to isolate the changed channels. Also, the plot of the real data shows that outliers appear occasionally in a few stocks. Therefore, it is better to combine with the Shewhart chart to detect the outliers or truncate off the outliers in order to detect the change segment that interests us the most. In addition, the plots show that the volatility change may be more important than the mean change, that will be left for future considerations. 

\section{Conclusion}

In this communication, we evaluated the performances of several common monitoring charts from the detection of transient signal point of view. When the length of the signal is unknown, the comparison of POD shows that EWMA and GRLT charts should be recommended. More importantly, as the EWMA chart does not require the reference values for the signal length, it is much more convenient to apply than other charts, particularly when the signal only appears in a small portion of the channels. 

As we only considered the case of mean change, more general models should be considered such as the variance (volatility) change in monitoring the financial data and correlated data case including both serial and cross-sectional dependence. In addition, the estimation of the signal length or the change point should be considered to estimate the true changed mean after the changed channels are isolated following an alarm. A study under multi-parameter exponential family model will be presented in a future communication. 

\begin{figure}[ht]
\begin{center}
\includegraphics[width =\textwidth,height=8.5in]{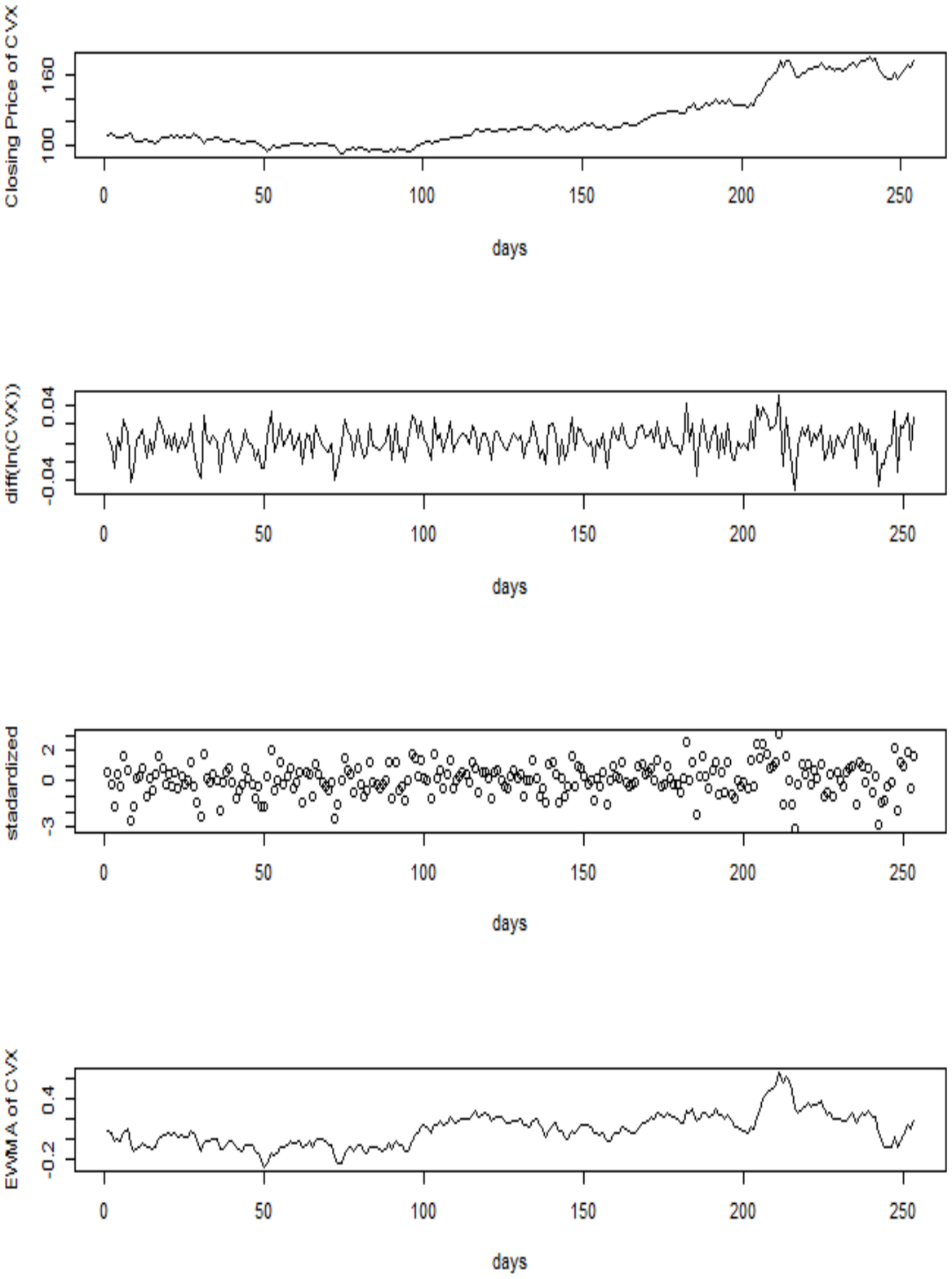}
Figure 2: EWMA chart for CVX
\end{center}
\end{figure}

\begin{figure}[ht]
\begin{center}
\includegraphics[width =\textwidth,height=8.5in]{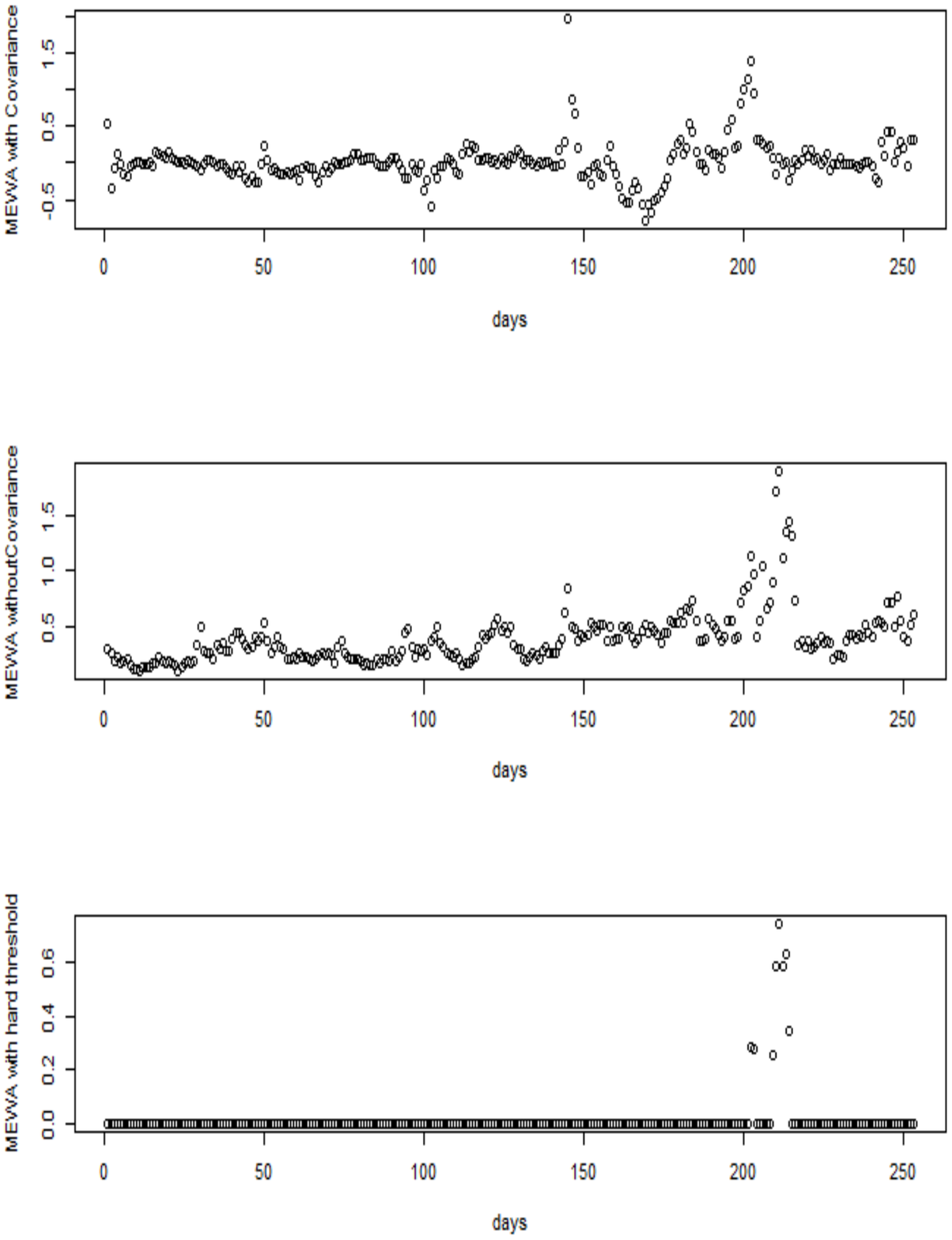}
Figure 3: MEWMA charts for 20 DJ30 stocks
\end{center}
\end{figure}

\section{Appendix}

\subsection {Proof of Theorem 1}

(i) Denote a changed measure $P_t^*(.)$ such that under $P_t^*(.)$, $X_0$ follows $N(0, \beta/(2-\beta))$, $X_t$ follows $N(\frac{\theta}{\beta}\sqrt{\frac{\beta}{2-\beta}}, 1)$, and $X_i$ follow $N(0,1) $ for $i\neq t $. Note that we can write
\[
\tau_{EW} =\inf \{t>0: Z_t= \frac{Y_t}{\sqrt{\beta/(2-\beta)}} >b\}.
\]
Then, the likelihood ratio of $P_t^*(.)$ with respect to $P_0^*(.)$ for $X_1,...X_L $ is
\[
\frac{dP_t^*}{dP_0^*}(Z_t) =\exp(\theta Z_t -\theta^2/2) ). 
\]
Select $\theta =b$. Denote by
\[
l_t =\theta Z_t -\frac{\theta^2}{2} =bZ_t -b^2/2,
\]
\[
\tilde{l}_t=\theta ( Z_t -\dot{\psi}(\theta))=b (Z_t -b).
\]
Then, we can write
\begin{eqnarray*}
P_0^*(\tau_{EW} \leq L) &=&P_0^* (\max_{1\leq t\leq L} Z_t \geq b  ) \\
&=& \sum_{t=1}^L E_0^*[ \frac{exp(l_t)}{\sum_{n=1}^L exp(l_n)}; \max_{1\leq t\leq L} Z_t \geq b ] \\
&=& \sum_{t=1}^L E_t^*[ \frac{1}{\sum_{n=1}^L exp(l_n)}; \max_{1\leq t\leq L} Z_t \geq b ] \\
&=& e^{-b^2/2} \sum_{t=1}^L E_t^*[ \frac{M_t}{S_t} e^{-(\tilde{l}_t +\max_{1\leq n \leq L} (l_n-l_t))}; \tilde{l}_t + \max_{1\leq n \leq L} (l_n-l_t))]\geq 0], 
\end{eqnarray*}
where $M_t =\max_{1\leq n \leq L} \exp(l_n-l_t))$ and $S_t =\sum_{1\leq n \leq L} \exp(l_n-l_t))$. 
Following the argument given in Siegmund and Yakir (2000), the two factors inside the expectations are approximately independent. Also, note that the variance of $\tilde{l}_t$ is approximately $b^2$.
\[
E_t^*[e^{-(\tilde{l}_t +\max_{1\leq n \leq L} (l_n-l_t))}; \tilde{l}_t +\max_{1\leq n \leq L} (l_n-l_t))\geq 0] \approx \frac{1}{b\sqrt{2\pi}}.
\]
Thus, when $L$ is large, 
\begin{equation}
P_0^*(\tau_{EW} \leq L) \approx \frac{L}{b\sqrt{2\pi}} e^{-b^2/2} \lim_{L \rightarrow \infty} E_t^*[ \frac{M_t}{S_t}]. 
\end{equation}
Locally, we can write 
\[
l_{t+1}-l_t = \theta( -\beta Z_t +\beta X_{t+1}/\sqrt{\beta/(2-\beta)} ) .
\]
As $\beta \rightarrow 0$, it has mean $\theta (-\beta \theta )=-\beta b^2 $ and variance approximately
$\theta^2 \beta^2 /(\beta/(2-\beta)) \approx 2\beta b^2$. Similarly, we can show the same results for $l_{t-1}-l_t$. That means $ l_n -l_t$ behaves approximately like a two-sided random walk as $\beta \rightarrow 0$ and 
\[
\lim_{L \rightarrow \infty} E_t^*[ \frac{M_t}{S_t}] \approx \beta b^2 \upsilon (b \sqrt{2\beta} ).
\]

(2) Here we define a changed measure $P_t^*(.)$ such that under $P_t^*(.)$, $X_{t-w+1} ,..., X_t$ are i.i.d following $N(\theta /\sqrt{w}, 1)$ and others are the same. Let $Z_t =(X_{t-w+1}+ ...+ X_t)/\sqrt{w}$ and $b=h\sqrt{w}$. Then the likelihood ratio of $X_1,..., X_L$ under $P_t^*(.)$ with respect to $P_0^*(.)$ equals to 
\[
\frac{dP_t^*}{dP_0^*}(Z_t) =exp(\theta Z_t -\theta^2/2 ).
\]
By taking $\theta =b $, we define
\[
l_t =\theta Z_t -\theta^2/2  = bZ_t -b^2/2,
\]
\[
\tilde{l}_t=\theta (Z_t -\theta)=b (Z_t -b). 
\]
Similar to (i), we can get the same form of approximation as in (2.3).  Note that we can write for $t-w\leq n \leq t$
\[
l_n -l_t =\frac{b}{\sqrt{w}} (\sum_{n-w+1}^n X_i - \sum_{t-w+1}^t X_i)=\frac{b}{\sqrt{w}}(\sum_{n-w+1}^{t-w}X_i -\sum_{n+1}^t X_i )
\]
It can seen that as $w \rightarrow \infty$, for $t-w\leq n \leq t+w$, $l_n-l_t$ is approximately a two-sided independent random walk with drift $ -\frac{b}{\sqrt{w}} \frac{\theta}{\sqrt{w}} = -\frac{b^2}{w}$ and variance $\frac{2b^2}{w}$ (as there are two random walks on each side one has drift 0 with variance $ \frac{b^2}{w}$ and another has drift $-\frac{b^2}{w}$ and variance $ \frac{b^2}{w}$). From Siegmund and Yakir (2000), we can approximate 
\[
\lim_{L \rightarrow \infty} E_t^*[ \frac{M_t}{S_t}] \approx \frac{ b^2}{w} \upsilon (b \sqrt{2/w} ).
\]
The expected result follows after some simplifications.   \qed

\subsection{Proof of Theorem 7}

We first define a changed measure $P_t^w(.)$ such that under the measure, $X_{t-w+1},..., X_t$ follows $N(\delta, I_N)$ and other $X_i$'s still follow $N(0, I_N)$. Then the log-likelihood ratio for $X_1,..., X_L$ under $P_t^w(.)$ with respect to $P_0^*(.)$ can be written as 
\[
l_t^w =-w \frac{||\delta||^2}{2} +\ln G( \frac{N}{2}, \frac{||\delta||^2 w^2 (\chi_t^w)^2}{4}),
\]
where $(\chi_t^w)^2 = \bar{X}_{t,w}^T\bar{X}_{t,w}$, and 
\[
G(a,y)=\sum_{i=0}^{\infty} \frac{\Gamma(a) y^i}{\Gamma(a+i) i!}
\]
see Basseville and Nikiforov (1993, Pg. 219). It is noted that from Abramowitz and Stegun (1964, Ch.13)
\[
G(a,y) = e^{-2\sqrt{y}} M(a-1/2, 2a-1, 4\sqrt{y}),
\]
where $M(a,b,z)$ is the confluent hypergeometric function and as $z\rightarrow \infty$, 
\[
M(a,b,z) =\frac{\Gamma(b)}{\Gamma(a)} e^z z^{a-b} (1+O(1/z)). 
\]
Thus,
\[
 \ln G( \frac{N}{2}, \frac{||\delta||^2 w^2 (\chi_t^w)^2}{4})
 \approx w ||\delta|| \chi_t^w -\frac{N-1}{2} \ln (2 w ||\delta|| \chi_t^w) +\ln \frac{\Gamma(N-1)}{\Gamma((N-1)/2)}.
 \]
Define the following stopping time
\[
\tau= \inf\{ t>0: \max_{w_0<w<w_1} l_t^w > ||\delta|| d +c \},
\]
where $c=\ln \frac{\Gamma(N-1)}{\Gamma((N-1)/2)}$. 
By using the change of measure technique, we can get
\[
P_0^*(\tau \leq L)= e^{-||\delta|| d-c }  \sum_{t=1}^L \sum_{w=w_0}^{w_1} E_t^w[ \frac{M_t}{S_t} e^{-(l_t^w-||\delta|| d -c+\ln M_t)} ; l_t^w -||\delta|| d -c + \ln M_t ]\geq 0],
 \]
where
$
M_t =\max_{1\leq u\leq L}\max_{w_0 <u<w_1} e^{(l_s^u-l_t^w)}$, and $ S_t = \sum_{1\leq u\leq L}\sum_{w_0 <u<w_1} e^{(l_s^u-l_t^w)}$. 

The key point is that the above approximation mainly relies on the w's that make the log-likelihood close to the boundary $|||\delta||d$. For large $w \chi_t^w$ with finite $N$,
\[
 l_t^w -c \approx w ||\delta|| \chi_t^w -w \frac{||\delta||^2}{2}  -\frac{N-1}{2} \ln (2 w ||\delta|| \chi_t^w)
 \]
Therefore, the CUSUM procedure defined by $\tau$ is equivalent to $\tau_{MCU}$. By the localization technique,  we have the approximation
\[
P_0^*(\tau_{MCU} \leq L) \approx L e^{-||\delta|| d -c} E_t^w[ \frac{M_t}{S_t}] \sum_{w=w_0}^{w_1} E_t^w[e^{-(l_t^w-||\delta|| d-c +\ln M_t)} ; l_t^w -||\delta|| b -c+ \ln M_t ]\geq 0]]. 
\]
Following the argument in Theorem 4, we only gives the main steps. First,  for large $w$, under $P_t^w(.)$, we can write
\begin{eqnarray*}
\chi_t^w &= &(\bar{X}_{t;w}^T\bar{X}_{t;w})^{1/2} =( (\delta +Z_{t;w}/\sqrt{w})^T (\delta +Z_{t;w}/\sqrt{w}))^{1/2} \\
&=& (\delta^T\delta +2\delta^TZ_{t;w}/\sqrt{w} + Z_{t;w}^TZ_{t;w}/w)^{1/2} \\
&=& ||\delta|| ( 1+ \frac{\gamma^TZ_{t;w}}{ ||\delta|| \sqrt{w}} + \frac{Z_{t;w}^TZ_{t;w}}{2w ||\delta||^2} -\frac{(\gamma^TZ_{t;w})^2}{2 w||\delta||^2} +o_p(1/w)) \\
&=& ||\delta|| +\frac{\gamma^TZ_{t;w}}{\sqrt{w}} + \frac{Z_{t;w}^TZ_{t;w}}{2w ||\delta||} -
\frac{(\gamma^TZ_{t;w})^2}{2 w||\delta||} +o_p(1/w)) ,
\end{eqnarray*}
where $\gamma =\delta/ ||\delta|| $ and $Z_{t;w} =\sqrt{w} \bar{X}_{t,w}$ that follows $N(0, I)$. From this, we have $E_t^w[\chi_t^w ] \approx ||\delta|| +\frac{N-1}{2w||\delta||}$ and 
$Var_t^w (\chi_t^w ) \approx 1/w$. 

From $w ||\delta||^2/2= ||\delta|| d$, we get $w^* = 2d/||\delta|| $. Thus, 
\[
E_t^w(l_t^w)\approx w ||\delta||^2/2 -\frac{N-1}{2} \ln (2w||\delta||^2) =w ||\delta||^2/2 - \frac{N-1}{2} \ln (4||\delta||d) . \]

The second term will contribute to the coefficient as a factor $ (4||\delta||d)^{(N-1)/2}$.  For general $w =w^*+\frac{2\mu}{||\delta||^2} \sqrt{w^*}$, we have
\[
l_t^w -||\delta||d (w^*+\frac{2\mu}{||\delta||^2} \sqrt{w^*})(||\delta|| (||\delta||+ \frac{\gamma^TZ_{t;w}}{\sqrt{w}} ) -||\delta||^2/2) -d||\delta|| \]
\[
\approx \mu \sqrt{w} +\sqrt{w} \gamma^TZ_{t;w}.
\]
Thus, 
\[
(l_t^w -||\delta|| d)/\sqrt{w} \rightarrow N(\mu, ||\delta||^2 ).
\]
The checking for the two-sided random walk behavior for $l_s^u-l_t^w$ needs a little more careful treatment. Under $P_t^w(.)$, for $ s<t, u>w$, we have
\begin{eqnarray*}
\chi_s^u &=& (\bar{X}_{s;u}^T\bar{X}_{s;u})^{1/2} =( (\frac{s-t+w}{u}\delta +Z_{s;u}/\sqrt{u})^T (\frac{s-t+w}{u}\delta +Z_{s;u}/\sqrt{u}))^{1/2} \\
&=& ( \frac{(s-t+w)^2}{u^2} ||\delta||^2 +2 \frac{s-t+w}{u} \frac{\delta^TZ_{s;u}}{\sqrt{u}} + \frac{Z_{s;u}^TZ_{s;u}}{u} )^{1/2} \\
&=& \frac{s-t+w}{u}  ||\delta|| +\frac{\delta^TZ_{s;u}}{\sqrt{u}}  +o_p(1/\sqrt{u}). 
\end{eqnarray*}
Thus, 
\[
l_s^u-l_t^w \approx u||\delta|| (\frac{s-t+w}{u} ||\delta|| + \frac{\delta^TZ_{s;u}}{\sqrt{u}} -||\delta||/2) -w||\delta|| (||\delta||/2 + \frac{\delta^TZ_{t;w}}{\sqrt{w}} )
\]
\[
=||\delta||^2(s-t+(w-u)/2) +\delta^T\sum_{i=s-u+1}^sZ_i -\delta^T\sum_{i=t-w+1}^tZ_i , 
\]
where $Z_i$ are i.i.d. $N(0, I)$ variables. Now we can see that it is approximately two independent two-sided random walk with drift $-||\delta||^2/2 $ and variance $||\delta||^2$. The result can be obtained as in the proof for Theorem 4. 
\qed

\section*{Acknowledgement} This work was done when the first author was visiting Department of Statistics, Stanford University. 

\section*{References}
\begin{description}
\item[{\rm Bauer, P. and Hackel, P. (1978)}] The use of MOSUMS for quality control. {\em Technometrics} {\bf 20(4)}: 431-436.
\item[{\rm Basseville, M., and Nikiforov, I. V. (1993)}] {\em Detection of Abrupt Changes: Theory and Applications}.  Prentice-Hall, Inc.    
\item[{\rm Gueppie, B. K., Fillatre, L, and Nikiforov, I. V. (2012)}]  Sequential detection of transient changes. {\em Sequential Analysis} {\bf 31(4)}: 528–547. 
\item[{\rm Knoth, S. (2021)}] Steady-state average run length(s) - methodology, formulas and numerics. {\em Sequential Analysis} {\bf 40(3)}: 405-426.
\item[{\rm Lai, T. L. (1974)}] Control charts based on weighted sums. {\em Annals of Statistics } {\bf 2}:head134-147. 
\item[{\rm Lai, T. L. (1995)}] Sequential changepoint detection in quality control and dynamical systems (with discussions). {\em Journal of Royal Statistical Society (B)} {\bf 57(4)}: 613-658.
\item[{\rm Lowry, C. A., Woodall, W. H., Champ, C. W., and Rigdon, S. E. (1992)}] A multivariate exponentially weighted moving average control chart. {\em Technometrics} {\bf34(1)}: 46-53.
\item [{\rm  Mei, Y. (2010)}] Efficient scalable schemes for monitoring a large number of data streams. {\em Biometrika} \textbf{97}: 419 -433.
\item[{\rm Ngai, H. M. and Zhang, J. (2001)}] Multivariate cumulative sum control charts based on projection pursuit. {\em Statistica Sinica} {\bf 11}: 747-766. 
\item[{\rm Noonan, J. and Zhigljavsky, A. (2020)}] Power of MOSUM test for online detection of transient change in mean. {\em Sequential Analysis} {\bf 39(2)}: 269-293.
\item[{\rm Page, E. S. (1954).}] Continuous inspection schemes. {\em Biometrika} {\bf 41}: 100-114.
 \item[{\rm  Pignatiello, J. J. and G. C. Runger, G. C. (1990)}]  Comparison of multivariate CUSUM charts {\em Journal of Quality Technology} {\bf 22:} 173-186.
 \item[{\rm Pollak, M. and Krieger, A. M. (2013)}] Shewhart revisited. {\em Sequential Analysis} {\bf 32}: 230–242. 
 \item[{\rm Reynolds, M. R., Jr. and Stoumbos, Z. G. (2004)}] Control charts and the efficient allocation of sampling resources. {\em Technometrics} {\bf 46(2)}: 200-214.
\item[{\rm Roberts, S. W. (1959)}] Control charts based on geometric moving average. {\em Technometrics} {\bf 1}:  239-250.
\item[{\rm Roberts, S. W. (1966)}] A comparison of some control chart procedures. {\em Technometrics} {\bf 8}: 411-430. 
\item[{\rm Shewhart, W. (1931)}] {\em Economic Control of Quality of Manufactured Products.} Princeton: Van Nostrand. 
\item[{\rm Shiryayev, A. N. (1963)}] On optimum methods in  quickest detection problems. {\em Theory of Probability and Applications} {\bf 13}: 22-46.
\item [{\rm Siegmund, D. (1985)}] {\it Sequential Analysis: Tests and Confidence Intervals}. Springer, New York.
\item [{\rm Siegmund, D. (1988)}] Approximate tail probabilities for the maxima of some random fields. {\em Annals of Probability} {\bf 16(2):} 487-501.
\item [{\rm Siegmund, D. and Venkatraman, E. S. (1995)}] Using the generalized likelihood ratio statistic for sequential detection of a change-point. {\em Annals of Statistics} {\bf 23(1):} 255-271.
\item[{\rm Siegmund, D. and Yakir, B.(2000)}]
Tail probabilities for the null distribution of scanning statistics. {\em Bernoulli } {\bf 6(2):} 191-213.
\item[{\rm Siegmund, D. and Yakir, B.(2008)}] Detecting the emergence of a signal in a noisy image. {\em Statistics and its Interface} {\bf 1}:3-12. 
\item[{\rm Siegmund, D., Yakir, B., and Zhang, N. (2010)}] Tail approximations for maxima of random fields by likelihood ratio transformations. {\em Sequential Analysis} {\bf 29}: 245–262.
\item[{\rm Siegmund, D., Yakir, B., and Zhang, N. R. (2011)}]
Detecting simultaneous variant intervals in aligned sequences. {\em Annals of Applied Statistics} {\bf 5(2A):} 645-668. 
\item[{\rm Srivastava, M. S. and Wu, Y. (1993)}] Comparison of CUSUM, EWMA, and Shiryayev-Roberts procedures for detection of a shift in the mean. {\em Annals of Statistics} {\bf 21(2)}: 645-670. 
\item [{\rm Tartakovsky, A. G. and Veeravalli, V. V. (2008)}] Asymptotically optimal quickest detection change detection in distributed sensor. {\em Sequential Analysis} \textbf{ 27}: 441-475.
\item[{\rm Wu, Y., and Wu, W. B. (2021).}] Sequential Detection of Common Transient Signals in High Dimensional Data Stream. To appear on {\em Naval Logistic Research Quarterly}
\item [{\rm Xie, Y. and Siegmund, D. (2013)}] Sequential multi-sensor change-point detection. {\em Annals of Statistics} \textbf{41}: 670-692.
\item[{\rm Xie, L., Xie, Y., and Moutstakides, G. V. (2021)}] Sequential subspace change-point detection. {\em Sequential Analysis} {\bf 39(3):} 307–335.
\end{description}

\end{document}